\documentclass[journal]{IEEEtran}
\usepackage{array}
\usepackage{marvosym}
\usepackage[cmex10]{amsmath}
\usepackage{amssymb,amsthm,epsfig}
\usepackage{graphpap}
\usepackage[latin1]{inputenc}
\usepackage{rotating}
\usepackage{algorithm}
\usepackage{setspace}
\usepackage{blkarray}
\usepackage{subfigure}
\usepackage{environ}
\usepackage{multirow}
\usepackage{setspace}
\usepackage{textcomp}
\usepackage[flushleft]{threeparttable}
\usepackage{color}
\usepackage{soul}

\usepackage{epsfig}
\usepackage{amsmath}
\usepackage{algcompatible}
\usepackage{kbordermatrix}
\usepackage{latexsym}
\usepackage{amsfonts}
\usepackage[latin1]{inputenc}
\usepackage{rotating}
\usepackage{multirow}
\usepackage{setspace}
\usepackage{textcomp}
\usepackage[pdftex]{hyperref}
\usepackage{bm}
\hypersetup{colorlinks=true, allcolors=black, linkcolor=black, citecolor=black }
\usepackage{mathtools}

\usepackage{mdwlist}
\usepackage{graphicx}
\usepackage{cite}
\usepackage{enumitem}

\usepackage{footnote}
\makesavenoteenv{tabular}
\usepackage{tabulary}
\usepackage[table]{xcolor}
\newcolumntype{P}[1]{>{\centering\arraybackslash}p{#1}}
\newcolumntype{M}[1]{>{\centering\arraybackslash}m{#1}}
\DeclareMathOperator*{\argmin}{arg\,min}
\begin{document}
\title{Repair and Resource Scheduling in Unbalanced Distribution Systems using Neighborhood Search}
%
\author{Anmar Arif,~\IEEEmembership{Student Member,~IEEE,}
       Zhaoyu Wang,~\IEEEmembership{Member,~IEEE,}
       Chen Chen,~\IEEEmembership{Member,~IEEE,}\\
       Jianhui Wang,~\IEEEmembership{Senior Member,~IEEE,}
\thanks{This project is sponsored in part by the U.S. Department of Energy Office of Electricity Delivery and Energy Reliability, the Iowa Energy Center, Iowa Economic Development Authority and its utility partners. The work is also supported by the U.S. Department of Energy's Solar Energy Technologies Office under Grant CPS\#34228.}
\thanks{A. Arif and Z. Wang are with the Department of Electrical and Computer Engineering, Iowa State University, Ames, IA, 50011. (Email:aiarif\MVAt iastate.edu,wzy\MVAt iastate.edu).}
\thanks{C. Chen is with the Energy Systems Division, Argonne National Laboratory, Lemont, IL 60439 USA (Email: morningchen\MVAt anl.gov)}
\thanks{J. Wang is with the Department of Electrical Engineering at Southern Methodist University, Dallas, TX, 75205 USA (email: jianhui\MVAt smu.edu)}
}
\maketitle
\begin{abstract}
This paper proposes an optimization strategy to assist utility operators to recover power distribution systems after large outages. Specifically, a mixed-integer linear programming (MILP) model is developed for co-optimizing crews, resources, and network operations. The MILP model coordinates damage isolation, network reconfiguration, distributed generator re-dispatch, and crew/resource logistics. {\color{black}In addition, a framework for integrating different types of photovoltaic (PV) systems in the restoration process is developed.} We consider two different types of crews, namely, line crews for damage repair and tree crews for obstacle removal. We also model the repair resource logistic constraints. Furthermore, a new algorithm is developed for solving the distribution system repair and restoration problem (DSRRP). The algorithm starts by solving DSRRP using an assignment-based method, then a neighborhood search method is designed to iteratively improve the solution. The proposed method is validated on modified IEEE 123- and {\color{black}8500-bus distribution test systems}. 
\end{abstract}

\begin{IEEEkeywords}
Outage management, power distribution system, repair crews, routing, service restoration
\end{IEEEkeywords}
\vspace{-0.2cm}
\section*{Nomenclature}
{
\footnotesize
\addcontentsline{toc}{section}{Nomenclature}
\begin{description}[style=multiline,leftmargin=3cm]
\item[\textbf{Sets and Indices}]
\end{description}
\begin{description}[style=multiline,leftmargin=1.5cm,itemsep=0.05cm] 
\item[$m/n$] Indices for damaged components and depots
\item[$c,r,w$] Index for crews, resources and depots
\item[$i/j$] Indices for buses
\item[$k$] Index for distribution line connecting $i$ and $j$
\item[$t,\varphi$] Index for time and phase number
\item[$C^L,C^T$] Set of line and tree crews
\item[$N$] Set of damaged components and the depot
\item[$N(c)$] Set of components assigned to crew $c$
\item[$\Omega_{B},\Omega_{P}$] Set of buses and depots
\item[$\Omega_{DK},\Omega_{DT}$] Set of damaged lines and lines damaged by trees.
{\color{black}\item[$\Omega_{ES},\Omega_{PV}$] Set of BESSs and PVs}
\item[$\Omega_{G},\Omega_{Sub}$] Set of buses with dispatchable generators and substations
\item[$\Omega_{K(.,i)}$] Set of lines with bus $i$ as the to bus
\item[$\Omega_{K(i,.)}$] Set of lines with bus $i$ as the from bus
\item[$\Omega_{K(l)}$] Set of lines in loop $l$
\item[$\Omega_{SW}$] Set of lines with switches
\end{description}
\begin{description}[style=multiline,leftmargin=3cm]
\item[\textbf{Parameters}]
\end{description}
\begin{description}[style=multiline,leftmargin=1.5cm,itemsep=0.06cm]
\item[$Cap^R_r$] The capacity required to carry resource $r$
\item[$Cap^C_c$] The maximum capacity of crew $c$
{\color{black}\item[$\underline{E}/\overline{E}^S_{i}$] The minimum/maximum energy state of BESS $i$
\item[$Ir_{i,t}$] Solar irradiance at bus $i$ and time $t$}
\item[$\mathcal{R}_{m,r}$] The number of type $r$ resources required to repair damaged component $m$
\item[$Res^D_{w,r}$] The number of type $r$ resources that are located in depot $w$
\item[$\rho^D_i,\rho^{SW}$] The cost of shedding the load at bus $i$ and cost of switching
\item[$M$] Large positive number
\item[$P/Q^D_{i,\varphi,t}$] Diversified active/reactive demand at bus $i$, phase $\varphi$ and time $t$
\item[$P/Q^U_{i,\varphi,t}$] Undiversified active/reactive demand at bus $i$ and phase $\varphi$
{\color{black}\item[$S,\bar{P}_i^{PV}$] The kVA and kW rating of PV $i$
\item[$S^{ES}_i$] The kVA rating of BESS $i$}
\item[$\mathcal{T}_{m,c}$] The estimated time needed to repair ( clear the trees at) damaged component $m$ for line (tree) crew $c$
\item[$tr_{m,n}$] Travel time between $m$ and $n$
\item[$\phi^0_c/\phi^1_c$] Start/End location of crew $c$
\item[$Z_k$] The impedance matrix of line $k$ 
\item[$\bm{p}_{k}$] Vector with binary entries for representing the phases of line $k$ 
\item[$\bm{a}_{k}$] Vector representing the ratio between the primary and secondary voltages for each phase of the voltage regulator on line $k$ 
\item[$\delta_{w,c}$] Binary parameter equals 1 if crew $c$ is positioned in depot $w$
{\color{black}\item[$\eta_c,\eta_d,\Delta t$] Charging and discharging efficiency, and the time step duration}
\end{description}
\begin{description}[style=multiline,leftmargin=3cm]
\item[\textbf{Decision Variables}]
\end{description}
\begin{description}[style=multiline,leftmargin=1.5cm,itemsep=0.08cm]
\item[$A^{L/T}_{m,c}$]Binary variable equal to 1 if component $m$ is assigned to line/tree crew $c$
\item[$Res^C_{c,w,r}$] Number of type $r$ resources that crew $c$ obtains from depot $w$
\item[$\gamma_{k,t}$] Binary variable indicates whether switch $k$ is operated in time $t$
\item[$\bm{S}_{k}$] A vector representing the apparent power of each phase for line $k$ at time $t$
\item[$\bm{U}_{i,t}$] A vector representing the squared voltage magnitude of each phase for bus $i$ at time $t$
\item[$\mathcal{X}_{i,t}$] Binary variable equal to 0 if bus $i$ is in an outage area at time $t$
\item[$E_{c,m,r}$]The number of type $r$ resources that crew $c$ has before repairing damaged component $m$
{\color{black}\item[$E^S_{i,t}$] Energy state of BESS $i$ at time $t$}
\item[$\alpha_{m,c}$] Arrival time of crew $c$ at damaged component $m$
\item[$f_{m,t}$] Binary variable equal to 1 if damaged component $m$ is repaired at time $t$
{\color{black}\item[$\mathcal{L}^{L},\mathcal{L}^{T}$] The expected times of the last repair conducted by the line and tree crews
\item[$P^{ch/dch}_{i,\varphi,t}$] Active power charge/discharge of the BESS at bus $i$}
\item[$P/Q_{i,\varphi,t}^L$] Active/reactive load supplied at bus $i$, phase $\varphi$ and time $t$
{\color{black}\item[$P/Q^{PV}_{i,\varphi,t}$] The active/reactive power output of the PV at bus $i$}
\item[$P/Q_{i,\varphi,t}^G$] Active/reactive power generated by DG at bus $i$, phase $\varphi$ and time $t$
\item[$P/Q^K_{k,\varphi,t}$] Active/reactive power flowing on line $k$, phase $\varphi$ and time $t$ 
\item[$\mathcal{P}_{c,w}$] A positive penalty term for the excess capacity that crew $c$ requires from depot $w$
{\color{black}\item[$\bar{tr}$] Maximum travel time for the crews }
\item[$u_{k,t}$] Binary variables indicating the status of the line $k$ at time $t$
{\color{black}\item[$u^{ES}_{i,t}$] Binary variable equals 1 if the BESS is charging and 0 for discharging
\item[$v^S_{i,t},v^f_{k,t}$] Virtual power generated at bus $i$ and the virtual flow on line $k$}
\item[$x_{m,n,c}$] Binary variable indicating whether crew $c$ moves from damaged components $m$ to $n$.
\item[$y_{i,t}$] Connection status of the load at bus $i$ and time $t$
\item[$z_{w,c}$] Binary variable equal to 1 if crew $c$ require additional resources from depot $w$
\end{description}
}
\section{Introduction}

\IEEEPARstart{T}HE combination of an aging electrical grid and a dramatic increase in severe storms has resulted in increasing large-scale power outages. In 2016,  the average outage duration for customers ranged from 27 minutes in Nebraska to 6 hours in West Virginia, while 20 hours in South Carolina due to Hurricane Matthew \cite{Outages2016}. The year 2017 experienced 18 major weather events around the world. The 2017 outages that were caused by hurricanes Harvey, Irma, and Maria alone have cost the U.S. around \$202 billion \cite{Outages2017}. 
{\color{black}Currently, utilities schedule the repairs using a list of predefined restoration priorities based on previous experiences, and network operation and repair scheduling are split into two different processes. This kind of approach does not capture the interdependence nature of the crew routing and network operation problems. Some customers cannot be served until the damaged lines are repaired, and the switching operation can affect the priorities of the repairs. Utilities commonly rely on the experiences of the operators. Our aim is to provide utilities with a better distribution system restoration decision-making process for coordinating crew scheduling, resource logistics, and network operations.}

{\color{black}Earlier work on distribution system restoration focused on network reconfiguration. In \cite{Jabr2012}, a mixed-integer conic program and mixed-integer linear program (MILP) were developed for network reconfiguration with the objective of minimizing the losses. The developed model included a spanning tree approach to enforce radiality and incorporated distributed generators (DGs). A MILP model and the genetic algorithm were used in \cite{Romero2005} for distribution network reconfiguration. The authors used graph theory to model the distribution network. 
Reference \cite{Hafez2018} proposed a decentralized agent-based method for service restoration. The developed approach divided the distribution system into several zones, where each zone was represented by an agent. The role of each agent was to maintain radial topology and operation limits and to maximize the served loads.}

{\color{black}Recent studies investigated the use of microgrids for distribution system restoration. The operation of multiple microgrids, with defined boundaries, in coordination with the distribution system has been investigated in \cite{Z_Wang2016} and \cite{Arif_MG}. The papers used stochastic programming for distribution system restoration with high penetration of DGs, including photovoltaic (PV) systems and battery energy storage systems (BESS). A decentralized method for coordinating networked microgrids and the distribution system was presented in \cite{Gao2018}. The authors modeled the operation of each microgrid as a second-order cone program and the coordination between the entities was achieved using the alternating direction method of multipliers algorithm. Other studies proposed sectionalizing the distribution network into microgrids; i.e., microgrids with dynamic boundaries. The authors in \cite{Chen2016} presented a MILP for microgrid formation of radial distribution networks to restore critical loads after outages. In \cite{Z_Wang2015}, the authors developed a two-stage stochastic mixed-integer nonlinear program to sectionalize the distribution network into multiple self-supplied microgrids. The paper included dispatchable DGs, such as microturbines and BESS, and PV systems. PVs and BESS were also considered in \cite{Trakas2018} for load restoration after wildfires.}

Although distribution system restoration has been long studied, there exist few efforts on integrating repair scheduling with recovery operation in power distribution systems. A pre-hurricane crew mobilization mathematical model was presented in \cite{pre-hurricane} for transmission networks. The authors used stochastic optimization to determine the number of crews to be mobilized to the potential damage locations. Also, the authors proposed a post-hurricane MILP model to assign repair crews to damaged components without considering the travel times and repair sequence. In \cite{Xu2007}, the authors developed a stochastic program that assigns crews to substations in order to inspect and repair the damage, but the approach neglected crew routing. The authors in \cite{Hent2015} presented a two-stage approach to decouple the crew routing and power restoration models in transmission systems. A MILP is solved in the first stage to find the priority of the damaged lines, and the routing problem is solved in the second stage using Constraint Programming. In \cite{Arif2016}, we developed a MILP that combines the distribution network operation and crew routing problems. The model was solved using a cluster-first route-second approach. Also, we developed a stochastic mixed integer linear program (SMIP) in \cite{Arif2018a} to solve the same problem with uncertainty. The problem was decomposed into two subproblems and solved using parallel progressive hedging. 

Several critical factors have been neglected in the previous work on this topic. First, when scheduling the crews, one must consider the different types of crews. There are mainly two types of crews: 1) line crews who are responsible for the actual repair of grid components; and 2) tree crews who remove obstacles in the damage sites before the line crews start the repairing work. The mathematical model for optimizing the crew schedule must include both types of crews to obtain an applicable solution. In terms of distribution system operation, the previous work did not include isolation of the damaged lines, which is imperative as the crews cannot repair a downed line until the power is cut off. {\color{black}Also, the connectivity of PV systems during outages in related work \cite{Arif_MG,Z_Wang2015,Trakas2018} does not represent the current practice. Due to technical, safety and regulatory issues, most on-grid (grid-tied) PV systems are disconnected during an outage (this is known as anti-islanding protection) \cite{Coastal_Solar}. On-grid PVs are required by law to have inverters with anti-islanding function \cite{Wholesale_Solar}.}

In this paper, we improve our previous work in \cite{Arif2016} and \cite{Arif2018a} by considering the 3-phase operation of the distribution network and modeling fault isolation constraints, coordinating tree and line crews, and resource logistics in the distribution system repair and restoration problem (DSRRP). {\color{black}Furthermore, a new framework for modeling different types of PV systems is developed. There are three main types of PV systems that are considered: 1) On-grid system: this type of PV is disconnected during an outage; 2) Hybrid on/off-grid (PV with BESS): the PV system operates on-grid in normal conditions, and off-grid during an outage (serves local load only); 3) PV + BESS with grid forming capabilities \cite{Kenning2018}: this system can restore part of the network that is not damaged if the fault is isolated. The idea of the proposed approach is to use a virtual network in parallel with the actual distribution network, and develop a mathematical formulation based on graph theory to identify the energized buses and the connectivity status of the PVs.

The crew routing problem is equivalent to the vehicle routing problem (VRP). VRP is an NP-hard combinatorial optimization problem that has been studied for a long time and remains challenging \cite{Braekers2016_vrp}. Combining VRP with the operation of distribution systems will further increase the complexity, therefore, some researchers opted to decouple the two problems \cite{Hent2015}.} In this paper, a tri-stage algorithm is developed to solve the proposed co-optimization model. The algorithm starts by solving an assignment problem, where the crews are assigned to the damaged components based on the expected working hours, distances between the crews and the outage locations, and the capacity of the crews. In the second stage, the DSRRP is solved with the crews dispatched to the assigned components from the first stage. In the third stage, a neighborhood search approach \cite{Fleszar2009} is used to iteratively improve the routing decisions obtained from stage two. The algorithm is used in a dynamically changing environment to handle the uncertainty of the repair time and other parameters. 
{\color{black}The contributions of this paper are summarized in the following:
\begin{itemize} 
 \item For the recovery operation of distribution systems, a mathematical formulation is developed for fault isolation and service restoration. Moreover, a formulation based on graph theory is developed for modeling the connectivity of PV systems during an outage.
 \item For crew routing, we model the coordination of line and tree crews as well as resource pick up. Equipment is needed to repair the damaged lines, however, a crew can only carry a limited number of supplies. Therefore, the crews need to go back to the depots and pick up additional supplies.   
 \item A new hybrid algorithm that combines mathematical programming and the neighborhood search method is designed to solve the computationally difficult repair and restoration problem. The algorithm is tested on modified IEEE 123- and IEEE 8500-bus distribution systems.
\end{itemize}
}
The rest of the paper is organized as follows. 
Section \ref{chap:4} develops the DSRRP mathematical formulation and Section \ref{chap:5} presents the algorithm for solving the model. The simulation results are presented in Section \ref{chap:6} and Section \ref{chap:7} concludes this paper.

\section{Distribution Network Repair and Restoration}\label{chap:4}

{\color{black}During extreme events, the outage management system (OMS) receives real-time data of the condition of the network from field devices, customer calls, and smart meters. Using the collected data, the OMS can estimate the locations of the outages, and the operator will dispatch field assessors to identify and document the exact locations of the damage. The DSRRP model can be incorporated in the OMS, where the model is solved to obtain the repair and restoration solution. The crew schedule is sent to the work management system (WMS), which communicates the tasks to the crews. The restoration plan and operations are sent to the distribution management system (DMS) and the system operator to confidently control the switches and DGs.}

In this paper, we assume that the assessors have located the damaged lines, and estimated the repair time and required resources. 
This section presents the mathematical model for coordinating line and tree crews, and the recovery operation of the network.
\subsection{Objective}
\begin{equation}
\small
\textrm{min}~  \sum \limits_{\forall t} \big {(} \mathop \sum_{\forall \varphi} \mathop \sum \limits_{\forall i}(1-y_{i,t})  \rho^D_i P^D_{i,\varphi,t}+ \rho^{SW}\sum_{\mathclap{k \in \Omega_{SW}}} \gamma_{k,t}\big {)}
\label{DSRRP_Obj}
\end{equation}
The first term in objective (\ref{DSRRP_Obj}) minimizes the cost of load shedding, while the second term minimizes the cost of operating the switches. The base load shedding cost is assumed to be \$14/kWh in this paper \cite{Ma2017}, and the base cost is multiplied by the load priority to obtain $\rho^D_i$. The switch operation cost is set to be \$8/time \cite{Malakar2013}.  
\subsection{Cold load pickup}
{
\small
\begin{equation}
P_{i,\varphi,t}^L=y_{i,t}P^D_{i,\varphi,t}+(y_{i,t}-y_{i,\rm{max}(t-\lambda,0)})P^U_{i,\varphi,t},~\forall i,\varphi,t
\label{P_clpu}
\end{equation}
\begin{equation}
Q_{i,\varphi,t}^L=y_{i,t}Q^D_{i,\varphi,t}+(y_{i,t}-y_{i,\rm{max}(t-\lambda,0)})Q^U_{i,\varphi,t},~\forall i,\varphi,t
\label{Q_clpu}
\end{equation}
\begin{equation}
{y_{i,t+1}} \ge {y_{i,t}}\;,\;\forall i,t
\label{units on}
\end{equation}
}
Constraints (\ref{P_clpu})-(\ref{Q_clpu}) set up the cold load pickup (CLPU) constraint \cite{Arif2018a}. In this paper, we employ two blocks to represent CLPU as suggested in \cite{Liu_2009}. The first block is for the undiversified load $P^U$ and the second for the diversified load $P^D$ (i.e., the steady-state load consumption). The use of two blocks decreases the computational burden imposed by nonlinear characteristics of CLPU and provides a conservative operation assumption to guarantee supply-load balance. Define $\lambda$ as the number of time steps required for the load to return to normal condition. The value of $\lambda$ is equal to the CLPU duration divided by the time step. The function max($t-\lambda,0$), is used to avoid negative time steps. If at time step $t=t_1$, a load goes from a de-energized state ($y_{i,t_1-1}=0$) to an energized one ($y_{i,t_1}=1$), it returns to normal condition at time step $t=t_1+\lambda$. $P_{i,\varphi,t}^U$ is added to $P_{i,\varphi,t}^D$ before time step $t_1+\lambda$ to represent the undiversified load. We assume that the duration of the CLPU decaying process is one hour \cite{Liu_2009}, and the total load at pickup time is $200\%$ of the steady state value \cite{Nagpal_2014}; i.e., $P_{i,\varphi,t}^U$ is set to be equal to $P_{i,\varphi,t}^D$. Constraint (\ref{units on}) indicates that once a load is served it cannot be shed.

\subsection{Power limits}
{\small
\begin{equation}
0 \le P_{i,\varphi,t}^{G} \le P_i^{G_{max}}\;,\;\forall i,\varphi,t
\label{DG limits}
\end{equation}
\begin{equation}
0 \le Q_{i,\varphi,t}^{G} \le Q_i^{G_{max}}\;,\;\forall i,\varphi,t
\label{QDG limits}
\end{equation}
\begin{equation}
- u_{k,t}P_k^{K_{min}} \le {P^K_{k,t}} \le u_{k,t}P_k^{K_{max}}\;,\;\forall k,t
\label{PL limits}
\end{equation}
\begin{equation}
- u_{k,t}Q_k^{K_{min}} \le {Q^K_{k,t}} \le u_{k,t}Q_k^{K_{max}}\;,\;\forall k,t
\label{QL limits}
\end{equation}
}
Constraints (\ref{DG limits})-(\ref{QL limits}) define the active and reactive power limits of the DGs and lines. The limits on the line-flow constraints are multiplied by $u_{k,t}$ so that if a line is damaged or a switch is opened, there will be no power flowing on it.

\subsection{Power flow equations}
{\small
\begin{equation}
\begin{split}
\mathop \sum_{\mathclap{\forall k \in K\left( {.,i} \right)}} {P^K_{k,\varphi,t}} + P^G_{i,\varphi,t} + &P^{PV}_{i,\varphi,t} + P^{dch}_{i,\varphi,t}=\\& \mathop \sum_{\mathclap{\forall k \in K\left( {i,.} \right)}} {P^K_{k,\varphi,t}} + P_{i,\varphi,t}^L + P^{ch}_{i,\varphi,t},\forall i,\varphi,t
\end{split}
\label{power balance}
\end{equation}
\begin{equation}
\begin{split}
\mathop \sum_{\mathclap{\forall k \in K\left( {.,i} \right)}} {Q^K_{k,\varphi,t}} + Q_{i,\varphi,t}^{G} + Q^{PV}_{i,\varphi,t} + &Q^{ES}_{i,\varphi,t} = \\& \mathop \sum_{\mathclap{\forall k \in K\left( {i,.} \right)}} {Q^K_{k,\varphi,t}} + Q_{i,\varphi,t}^L,\forall i,\varphi,t
\end{split}
\label{reactive balance}
\end{equation}
\begin{equation}
\begin{split}
{\bm{U}_{j,t}} - {\bm{U}_{i,t}} + \bm{\bar{Z}}_k &\bm{S}_{k}^*+  \bm{\bar{Z}}_k^*  \bm{S}_{k}\le \\ & (2-u_{k,t}-\bm{p}_{k}) M,\forall k \in \Omega_L,t
\end{split}
\label{voltage1}
\end{equation}
\begin{equation}
\begin{split}
{\bm{U}_{j,t}} - {\bm{U}_{i,t}} + &\bm{\bar{Z}}_k \bm{S}_{k}^*+  \bm{\bar{Z}}_k^*  \bm{S}_{k}\geq \\ & -(2-u_{k,t}-\bm{p}_{k}) M,\forall k \in \Omega_L,t
\end{split}
\label{voltage2}
\end{equation}
}
Constraints (\ref{power balance})-(\ref{reactive balance}) are 3-phase active and reactive power node balance constraints. Constraints (\ref{voltage1})-(\ref{voltage2}) represent Kirchhoff's voltage law. $\textbf{S}_{i,j} \in \mathbb{C}^{3 \times 1}$ is the three-phase apparent power from bus $i$ and $j$, and $\textbf{U}_i = [|V_i^a|^2,|V_i^b|^2,|V_i^c|^2]^T$. The matrix $\bar{Z}_{i,j}$ equals $\textbf{A} \odot \textbf{Z}_{i,j}$, where $\textbf{Z}_{i,j} \in \mathbb{C}^{3 \times 3}$ is the impedance matrix of the line, and $\textbf{A}$ is a phase shift matrix. Detailed derivation of (\ref{voltage1}) and (\ref{voltage2}) is provided in \cite{B_Chen2018}. The big $M$ method is used to decouple the voltages between lines that are disconnected or damaged. Also, if line $k(i,j)$ is two-phase (e.g., phases $a$ and $c$), then the voltage constraint is only applied to these two phases, which is realized by including $\bm{p_k}$. The vector $\bm{p_k} \in \{0,1\}^{3 \times 1}$ represents the phases of line $k$; e.g., for line $k$ with phases $a,c$, $\bm{p_k}=[1,0,1]$.
%

\subsection{Reconfiguration and Isolation}
{\small
\begin{equation}
\mathcal{X}_{i,t} U_{min}  \le {\bm{U}_{i,t}} \le \mathcal{X}_{i,t} U_{max} \;,\;\forall i,t
\label{V limits}
\end{equation}
\begin{equation}
2u_{k,t} \ge \mathcal{X}_{i,t} + \mathcal{X}_{j,t}, \forall k \in \Omega_{DK}, t
\label{xx=0}
\end{equation}
\begin{equation}
u_{k,t} = 1, \forall k \not \in \{\Omega_{SW} \cup \Omega_{DK}\},t
\label{uk=1}
\end{equation}
\begin{equation}
\sum_{k \in \Omega_{K(l)}} u_{k,t} \le |\Omega_{K(l)}|-1, \forall l,t
\label{radial_con}
\end{equation}
\begin{equation}
\gamma_{k,t} \ge u_{k,t}-u_{k,t-1}, \forall k \in \Omega_{SW}, t
\label{sw1}
\end{equation}
\begin{equation}
\gamma_{k,t} \ge u_{k,t-1}-u_{k,t}, \forall k \in \Omega_{SW}, t
\label{sw2}
\end{equation}
}
Constraint (\ref{V limits}) ensures that the voltage is within a specified limit, and is set to equal to 0 if the bus is in an on-outage area. Constraint (\ref{xx=0}) sets the values of $\mathcal{X}_i$ and $\mathcal{X}_j$ to be 0 if the line is damaged, therefore, the voltages on the buses between damaged lines are forced to be 0 using constraint (\ref{V limits}). Subsequently, the zero voltage propagates on the rest of the network through constraints (\ref{voltage1}) and (\ref{voltage2}) until a circuit breaker (CB) or sectionalizer stops the propagation. If the voltages on two connected buses are zero, then the power flow is forced to be zero through constraints (\ref{voltage1}) and (\ref{voltage2}). Constraint (\ref{uk=1}) defines the default status of the lines that are not damaged or not switchable. Constraint (\ref{radial_con}) is the radiality constraint. Radiality is enforced by introducing constraints for ensuring that at least one of the lines of each possible loop in the network is open \cite{Borghetti2012}. A depth-first search method is used to identify the possible loops in the network and the lines associated with them. Constraint (\ref{sw1})-(\ref{sw2}) are used in order to limit the number of switching operations. {\color{black}We assume that all switches are remotely controllable.} Let $\gamma_{k,t}$ equal to 1 if the line switches its status from 0 (off) to 1 (on), or 1 (on) to 0 (off). This variable is included in the objective to minimize the number of switching operations.  
{\color{black}
\subsection{PV Systems}
In this study, we consider three types of PV systems:
\begin{itemize}
\item Type 1: on-grid (grid-tied) PV ($\Omega_{PV}^{G}$): during an outage, the PV is switched off. This type of PV is the most commonly used one especially for residential customers \cite{Meehan2017}. The on-grid system uses a standard grid-tied inverter and does not have any battery storage.
\item Type 2: hybrid on-grid/off-grid PV + BESS ($\Omega_{PV}^{H}$): this system is an on-grid system that can disconnect from the grid after an outage and uses battery backup supply.
\item Type 3: grid-forming PV + BESS ($\Omega_{PV}^{C}$): this system is an on-grid system that can support a large section of the network \cite{Kenning2018}. After an outage, the PV and battery system can provide energy to the healthy parts of the network.
\end{itemize} 

\subsubsection{PV Active and Reactive Power}
The active and reactive powers of a PV depend on the rating of the solar cell and the solar irradiance. The active output power from the PVs is determined using constraints (\ref{PV_out}) and (\ref{PV_out2}). The PV inverters can provide reactive power support, which is constrained by (\ref{QPV_limit1}) and (\ref{QPV_limit2}) \cite{Q_Zhang2019}.
{\small
\begin{equation}
P^{PV}_{i,\varphi,t} = \frac{Ir_{i,t}}{(1000W/m^2)} \overline{P}^{PV}_i, \forall i \in \Omega_{PV}\backslash \Omega_{PV}^{G}, \varphi, t
\label{PV_out} 
\end{equation}
\begin{equation}
P^{PV}_{i,\varphi,t} = \mathcal{X}_{i,t} \frac{Ir_{i,t}}{(1000W/m^2)} \overline{P}^{PV}_i, \forall i \in \Omega_{PV}^{G}, \varphi, t
\label{PV_out2} 
\end{equation}
\begin{equation}
|Q^{PV}_{i,\varphi,t}| \le \sqrt{(S^{PV}_i)^2-(\hat{P}^{PV}_{i,t})^2}, \forall i \in \Omega_{PV}\backslash \Omega_{PV}^{G}, \varphi,t
\label{QPV_limit1}
\end{equation}
\begin{equation}
|Q^{PV}_{i,\varphi,t}| \le \mathcal{X}_{i,t}~\sqrt{(S^{PV}_i)^2-(\hat{P}^{PV}_{i,t})^2}, \forall i \in \Omega_{PV}^{G}, \varphi,t
\label{QPV_limit2}
\end{equation}
\begin{equation*}
\text{where} ~~~~~~ \hat{P}^{PV}_{i,t} = \frac{Ir_{i,t}}{(1000W/m^2)} \overline{P}^{PV}_i 
\end{equation*}
}
PVs of types $\Omega_{PV}^{H}$ and $\Omega_{PV}^{C}$ are able to disconnect from the grid and serve the on-site load. On the other hand, on-grid PVs are disconnected and the on-site load is not served by the PVs during an outage, therefore, the right-hand side in (\ref{PV_out2}) and (\ref{QPV_limit2}) are multiplied by $\mathcal{X}_{i}$. Note that $|f(x)| \le x$ is equivalent to $-x \le f(x) \le x$.
\subsubsection{PV Connectivity}
In this paper, we assume that the network can be restored using the grid-forming sources in $\Omega_{PV}^C \cup \Omega_G \cup \Omega_{Sub}$. A PV of type $\Omega_{PV}^{G}$ or $\Omega_{PV}^{H}$ can connect to the grid only after the PV bus is energized. Consider the network shown in Fig. \ref{Grid-forming}. Due to a line damage, the network is divided into four islands. Island A can be energized by the substation, therefore, the PV at bus 10 can be connected with the grid. Island B must be isolated because of the damaged line. Island C does not have any grid-forming generators; hence, it will not be active and the grid-tied PV will be disconnected. However, the PV+BESS system at bus 7 can energize the local load. Island D can be energized by the grid-forming PV+BESS system at bus 4.

\begin{figure}[htpb!]
\setlength{\abovecaptionskip}{0pt} 
\setlength{\belowcaptionskip}{0pt} 
  \centering
    \includegraphics[width=0.49\textwidth]{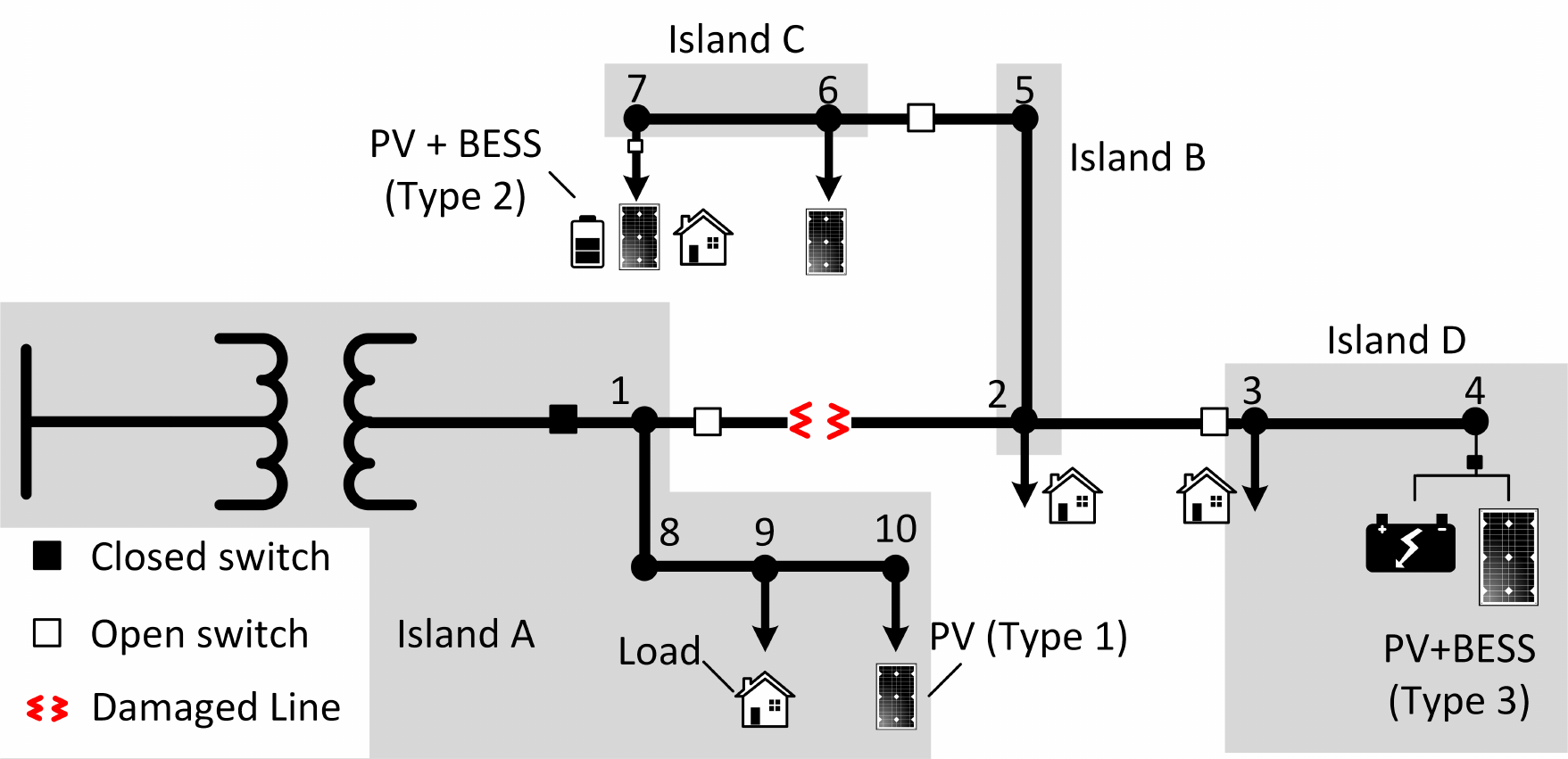}
  \caption{A single line diagram of a network with one damaged line.}\label{Grid-forming}
\end{figure}

The connectivity constraints of the PVs are represented by constraints (\ref{virtual_flow})-(\ref{xx>y_2}). The idea of the approach is to use virtual sources, loads, and flow to identify the energized buses in the network. The constraints for the virtual framework are formulated as follows:
{\small
\begin{equation}
v_{i,t}^S + \sum_{k \in K(.,i)}v^f_{k,t} = \mathcal{X}_{i,t}+ \sum_{k \in K(i,.)}v^f_{k,t}, \forall i, t
\label{virtual_flow}
\end{equation}
\begin{equation}
\sum_{\forall t} v^S_{i,t} = 0, \forall i \in \Omega_B \backslash \{\Omega_{PV}^C \cup \Omega_G \cup \Omega_{Sub}\}
\label{vs_limit}
\end{equation}
\begin{equation}
-u_{k,t}M \le v^f_{k,t} \le u_{k,t}M, \forall k \in \Omega_K, t
\label{vf_limit}
\end{equation}
\begin{equation}
\mathcal{X}_{i,t} \ge y_{i,t}, \forall i \in \Omega_B\backslash \{\Omega_G \cup \Omega_{PV}^C \cup \Omega_{PV}^H\}, t
\label{xx>y_2}
\end{equation}
}
To identify whether an island is energized by grid-forming generators or not, we create a virtual network. First, each grid-forming generator is replaced by a virtual source/generator with infinite capacity. Other power sources without grid-forming capability (e.g., grid-tied PVs) are removed. Also, virtual loads with magnitude of 1 are placed on each bus, and the actual loads are removed. For example, the network shown in Fig. \ref{Grid-forming} is transformed to the network shown in Fig. \ref{Grid-forming_V}. In the mathematical model, we add a node-balance equation for each virtual bus. If the virtual load at a bus is served, then that bus is energized. Therefore, for islands without grid-forming generators, all buses will be de-energized as the virtual loads in the island cannot be served. 
Constraint (\ref{virtual_flow}) is the node balance constraint for the virtual network. Constraints (\ref{vs_limit}) states that buses without grid-forming power generators do not have virtual sources. The variable $v^f_k$ represents the virtual flow on line $k$ and each bus is given a load of 1 that is multiplied by $\mathcal{X}_i$. Therefore, $\mathcal{X}_i$ = 1 (bus $i$ is energized) if the virtual load can be served by a virtual source and 0 (bus $i$ is de-energized) otherwise. The virtual flow limits are defined in (\ref{vf_limit}). If bus $i$ is de-energized, then the load must be shed (\ref{xx>y_2}), unless bus $i$ has a local power source.
\begin{figure}[htpb!]
\setlength{\abovecaptionskip}{0pt} 
\setlength{\belowcaptionskip}{0pt} 
  \centering
    \includegraphics[width=0.49\textwidth]{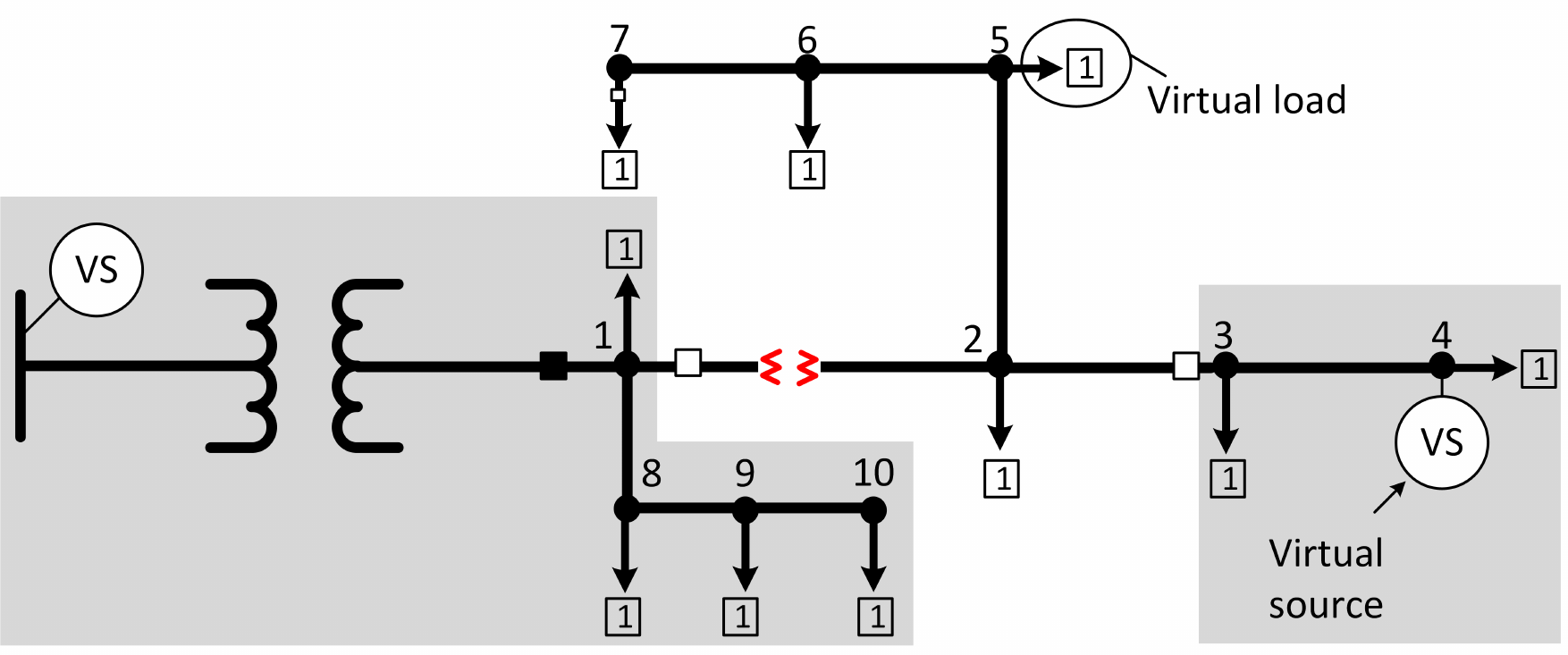}
  \caption{A virtual network created for the network shown in Fig. \ref{Grid-forming}.}\label{Grid-forming_V}
\end{figure}

\subsection{BESS}
{\small
\begin{equation}
0 \le {P^{ch}_{i,\varphi,t}} \le u^{ES}_{i,t}\overline{P}^{ch}_i,\;\forall i \in \Omega_{ES},\varphi,t
\label{Pch_limits}
\end{equation}
\begin{equation}
0 \le {P^{dch}_{i,\varphi,t}} \le (1-u^{ES}_{i,t})\overline{P}^{dch}_i,\forall i \in \Omega_{ES},\varphi,t
\label{Pdch_limits}
\end{equation}
\begin{equation}
{E}^{S}_{i,t} = {E}^{S}_{i,t-1}+\Delta t (\eta_c \sum_{\forall \varphi} P^{ch}_{i,\varphi,t}- \frac{\sum_{\forall \varphi} P^{dch}_{i,\varphi,t}}{\eta_d}), \forall i \in \Omega_{ES},t
\label{SOC_next}
\end{equation}
\begin{equation}
\underline{E}_{i}^{S} \le {E}_{i,t}^S \le \overline{E}_{i}^{S}, \forall i \in \Omega_{ES}, t
\label{SOC_limits}
\end{equation}
\begin{equation}
(Q^{ES}_{i,\varphi,t})^2+(P^{ch}_{i,\varphi,t}+P^{dch}_{i,\varphi,t})^2 \le ({S^{ES}_{i}})^2, \forall i \in \Omega_{ES}, \varphi,t
\label{Qcharging_limit1}
\end{equation}
}
Binary variable $u^{ES}$ represents the charging (1) and discharging (0) state of the BESS. Limits on the charge and discharge powers are imposed using constraints (\ref{Pch_limits}) and (\ref{Pdch_limits}), respectively. Constraint (\ref{SOC_next}) represents the dynamic state of energy for each BESS, where the efficiencies $\eta_c$ and $\eta_d$ are assumed to be 0.95. The energy is limited to a minimum and maximum value in (\ref{SOC_limits}). $E^S_{i,t}$ is assumed to be between 0.2 and 0.9 of the rated capacity in this paper. The active and reactive power should not exceed the rating of the BESS, as enforced by (\ref{Qcharging_limit1}) \cite{Abdeltawab2017}. Constraint (\ref{Qcharging_limit1}) is quadratic, therefore, it is linearized using the circular constraint linearization method presented in \cite{X_Chen2016}. Subsequently, constraint (\ref{Qcharging_limit1}) is replaced by (\ref{Qch_limits})-(\ref{Sch_linear2}).
{
\small
\begin{equation}
- S^{ES}_i \le {Q^{ES}_{i,\varphi,t}} \le S^{ES}_i,\forall i \in \Omega_{ES},\varphi,t
\tag{\ref{Qcharging_limit1}a}\label{Qch_limits}
\end{equation}
\begin{equation}
|{(P^{ch}_{i,\varphi,t}+P^{dch}_{i,\varphi,t}) + Q^{ES}_{i,\varphi,t}}| \le \sqrt{2} S^{ES}_i,\forall i \in \Omega_{ES},\varphi,t
\tag{\ref{Qcharging_limit1}b}\label{Sch_linear1}
\end{equation}
\begin{equation}
|{(P^{ch}_{i,\varphi,t}+P^{dch}_{i,\varphi,t}) - Q^{ES}_{i,\varphi,t}}| \le \sqrt{2}S^{ES}_i,\forall i \in \Omega_{ES},\varphi,t
\tag{\ref{Qcharging_limit1}c}\label{Sch_linear2}
\end{equation}
}
} 
\vspace{-0.7cm}
\subsection{Routing Constraints}

{\color{black}The routing problem can be defined by a complete graph with nodes and edges $G(N,E)$. The node set $N$ in the undirected graph contains the depot and damaged components, and the edge set ${E = \left\lbrace(m,n) | m,n \in N; m \neq n \right\rbrace}$ represents the edges connecting each two components. The graph $G$ can be obtained from a transportation network ($\hat{G}$). Transportation networks can be represented by nodes (depots, damaged components, intersection nodes) and paths connecting the nodes. Consider the transportation network shown in Fig. \ref{fig_transportation}a, where there are two damaged components and one depot. The information that is required by the DSRRP model is the travel time between the damaged components and the depot. Therefore, we can convert $\hat{G}$ to the network $G$ shown in Fig. \ref{fig_transportation}b by finding the shortest paths between damaged components and the depot \cite{Lei_2019}, which can be obtained using shortest path algorithms such as Dijkstra's algorithm \cite{Dijkstra}. In the example shown in Fig. \ref{fig_transportation}, the shortest path between the depot and damaged component A has a total length of 3 units. Therefore, the depot is connected directly to damaged component A in $G$ with a length of 3 units. The same procedure is conducted to form the rest of the network $G$. If a path between two nodes in $\hat{G}$ is completely blocked or severely damaged, then the travel time of the path can be set to a large value $|T|$, where $T$ is the time horizon. In practice, utilities use geographic information system (GIS) software to map the distribution network. Real-time data about road conditions, location of the crews, and status of the components are fed into the GIS. The utilities can then use the GIS to estimate the travel times.

\begin{figure}[h!]
\setlength{\abovecaptionskip}{0pt} 
\setlength{\belowcaptionskip}{0pt} 
  \centering
    \includegraphics[width=0.45\textwidth]{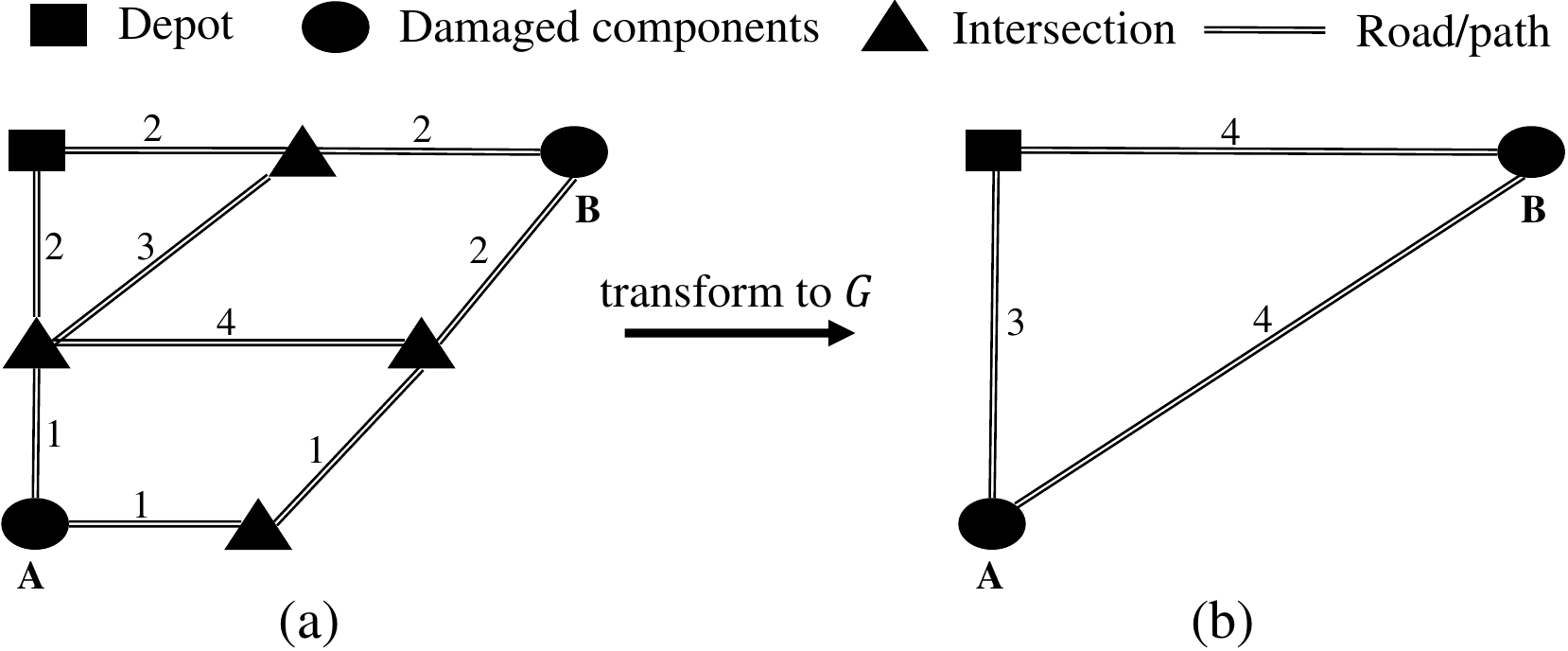}
  \caption{Example of (a) a transportation network transformed to (b) graph $G$ for the crew routing model.}\label{fig_transportation}
\end{figure} 

Our purpose is to find an optimal route for each crew to reach the damaged components. The value of $x_{m,n,c}$ determines whether the path crew $c$ travels includes the edge $(m,n)$ with $m$ preceding $n$. The routing constraints for the first stage problem are formulated as follows:}

{\small
\begin{equation}
\mathop \sum \limits_{\forall m \in N}{x_{\phi^0_c,m,c}} = 1, \forall c
\label{start_from_depot}
\end{equation}
\begin{equation}
\mathop \sum \limits_{\forall m \in N} {x_{m,\phi^1_c,c}} = 1, \forall c
\label{back_to_depot}
\end{equation}
\begin{equation}
\mathop \sum_{\mathclap{\forall n \in N\backslash \left\{ m \right\}}} {x_{m,n,c}} - \mathop \sum_{\mathclap{{\forall n \in N\backslash \left\{ m \right\}}}} {x_{n,m,c}} = 0,\forall c,m \in N\backslash \left\{\phi^0_c,\phi^1_c\right\}
\label{x=0}
\end{equation}
\begin{equation}
\mathop \sum \limits_{\forall c \in C^L} \mathop \sum \limits_{\forall m \in N\backslash \left\{n\right\}}x_{m,n,c}=1, \forall n \in \Omega_{DK}
\label{x=1_L}
\end{equation}
\begin{equation}
\mathop \sum \limits_{\forall c \in C^T} \mathop \sum \limits_{\forall m \in N\backslash \left\{n\right\}}x_{m,n,c}=1, \forall n \in \Omega_{DT}
\label{x=1_T}
\end{equation}}Constraint (\ref{start_from_depot})-(\ref{back_to_depot}) guarantee that each crew starts and ends its route at the defined start ($\phi^0_c$) and end ($\phi^1_c$) locations. Constraint (\ref{x=0}) is the flow conservation constraint; i.e., once a crew arrives at a damaged component, the crew moves to the next location after finishing the repairs. Constraint (\ref{x=1_L}) ensures that each damaged component is repaired by only one line crews, while (\ref{x=1_T}) ensures that each damaged component that needs removing a fallen tree first, is assigned to one tree crew.
\subsection{Arrival Time}
{\small
\begin{equation}
\begin{split}
\alpha_{m,c} + {\mathcal{T}_{m,c}} +& t{r_{m,n}} - \left( {1 - {x_{m,n,c}}} \right)M \le \alpha_{n,c}\\ &\forall m \in N \backslash \{\phi^1_c\},n \in N\backslash \left\{\phi^0_c,m\right\},c
\end {split}
\label{AT1}
\end{equation}
\begin{equation}
\sum_{c \in C^L}\alpha_{m,c} \ge \sum_{c \in C^T}\alpha_{m,c}+\mathcal{T}_{m,c}\sum \limits_{\mathclap{\forall n \in N}} x_{m,n,c}, \forall m \in \Omega_{DT}
\label{AT_L>AT_T}
\end{equation}
}
Constraint (\ref{AT1}) is used to calculate the arrival time (the time when crew $c$ starts repairing component $m$) for each crew at each damaged component. For a crew that travels from damaged component $m$ to $n$, $\alpha_{n,c}$ equals $\alpha_{m,c} + {\mathcal{T}_{m,c}} + tr_{m,n}$. Big $M$ is used to decouple the times to arrive at components $m$ and $n$ if the crew does not travel from $m$ to $n$. Constraint (\ref{AT_L>AT_T}) indicates that the line crews start repairing the damaged components after the tree crews clear the obstacles.
 
\subsection{Resource and Pick Up Constraints}
{\small
\begin{equation}
Res^D_{w,r} \ge ~~\sum_{\mathclap{\forall c \in C^L, \phi^0_c=w}}~Res^C_{c,\phi^0_c,r} + \sum_{\forall c \in C^L}Res^C_{c,w,r}, \forall w,r
\label{Res1}
\end{equation}
\begin{equation}
\sum_{\forall r}Cap^R_r E_{c,m,r} \le Cap^C_c, \forall m,c \in C^L
\label{Res2}
\end{equation}
\begin{equation}
\sum_{\forall n \in N}x_{n,m,c} \mathcal{R}_{m,r} \le E_{c,m,r}, \forall m,r,c \in C^L
\label{Res3}
\end{equation}
\begin{equation}
\begin{split}
&-M(1-x_{m,n,c}) \le \\
&E_{c,m,r}-\mathcal{R}_{m,r}-E_{c,n,r}\le M(1-x_{m,n,c}),\\
&\forall m \in N \backslash \{\phi^1_c\},n \in N\backslash \left\{\phi^0_c,m\right\},c \in C^L,r
\end {split}
\label{Res4}
\end{equation}
\begin{equation}
\begin{split}
&-M(1-x_{w,n,c})\le E_{c,w,r}+Res_{c,w,r}^C-E_{c,n,r} \le \\& M(1-x_{w,n,c}), \forall w, n \in N\backslash \left\{\phi^0_c,\phi^1_c,w\right\},c \in C^L,r
\end {split}
\label{Res5}
\end{equation}
\begin{equation}
\begin{split}
-M(1-&x_{\phi^0_c,n,c})\le Res_{c,\phi^0_c,r}^C-E_{c,n,r}\le\\& M(1-x_{\phi^0_c,n,c}), \forall n \in N\backslash \left\{\phi^0_c\right\},c \in C^L,r
\end {split}
\label{Res6}
\end{equation}
}
Constraint (\ref{Res1}) states that the total resources that the crews obtain from depot $w$ must be less or equal to the amount of available resources in the depot. The amount of resources that a crew can carry must be limited by the crew's capacity, which is realized by constraint (\ref{Res2}). Constraint (\ref{Res3}) indicates that the crews must have enough resources to repair the damaged components. Constraint (\ref{Res4}) ensures that if a crew travels from $m$ to $n$, then the resources that the crew have when arriving at location $n$ is $E_{c,n,r}=E_{c,m,r}-\mathcal{R}_{m,r}$. If a crew goes to depot $w$ to pick up supplies and travels to damaged component $n$, then $E_{c,n,r}=E_{c,w,r}+Res_{c,w,r}^C$, which is enforced by (\ref{Res5}). Constraint (\ref{Res6}) ensures that the number of resources that the crew has at the first damaged component is equal to the resources obtained at the starting location.
\subsection{Restoration Time}
{\small
\begin{equation}
\mathop \sum \limits_{\forall t} {f_{m,t}} = 1\;,\;\forall m \in \Omega_D
\label{f=1}
\end{equation}
\begin{equation}
\mathop \sum \limits_{\forall t} t{f_{m,t}} \ge \mathop \sum \limits_{\forall {c}} ( \alpha_{m,c} + {\mathcal{T}_{m,c}}\mathop \sum \limits_{\mathclap{\forall n \in N}} x_{m,n,c} ), \forall m \in \Omega_D
\label{tf1}
\end{equation}
\begin{equation}
0 \le \alpha_{m,c} \le M \mathop \sum \limits_{n \in N}x_{n,m,c},\;\forall m \in N\backslash \left\{ {\phi^0_c,\phi^1_c} \right\},c
\label{AT2}
\end{equation}
\vspace{-0.2cm}
\begin{equation}
{u_{m,t}} = \mathop \sum \limits_{\tau  = 1}^{t} {f_{m,\tau}}\;,\;\forall m \in \Omega_{DL},t
\label{z=f}
\end{equation}
\begin{equation}
\{f,x,u,y,\mathcal{X},\gamma\} \in \{0,1\}, \{E,Res^C\} \ge 0
\label{DSRRP_var_def}
\end{equation}
}
Constraints (\ref{f=1})-(\ref{z=f}) are used to connect the crew scheduling and power operation problems. Let $f_{m,t}$ denote the time when the damaged component is repaired by the line crews, which equals 1 in one time interval as enforced by (\ref{f=1}). Equation (\ref{tf1}) determines the time when a damaged component is repaired by setting $\sum_{\forall t} t{f_{m,t}}$ to be greater than or equal to $\alpha_{m,c}+\mathcal{T}_{m,c}$ of the crew assigned to damaged component $m$. Constraint (\ref{AT2}) is used to set $\alpha_{m,c}=0$ if crew $c$ does not travel to component $m$, so it would not affect constraint (\ref{tf1}). Finally, constraint (\ref{z=f}) indicates that the restored component becomes available after it is repaired, and remains available in all subsequent periods. For example, if $f_{m,t}=[0,0,1,0,0,0]$ then $u_{m,t}=[0,0,1,1,1,1]$. 

\section{Solution Algorithm}\label{chap:5}

A three-stage algorithm for solving the combined routing and distribution system operation problem is presented in this section, where the stages are: assignment, initial solution, and neighborhood search. 
Furthermore, to compare the developed method with current practices, a priority-based method that mimics the utilities' scheduling procedures is developed. 
\subsection{Reoptimization Algorithm} 
\subsubsection{Assignment}
By assigning the damaged components to the crews, the large VRP problem can be converted to multiple small-size Travelling Salesman Problems (TSP) \cite{laporte_2009}. The assignment problem is formulated as follows:
{
\small
{\color{black}
\begin{equation}
\min \mathcal{L}^L + \mathcal{L}^T + \sum_{\forall c}\sum_{\forall w}\mathcal{P}_{c,w} + \bar{tr}
\label{Pre_obj}
\end{equation}
\begin{equation}
\mathcal{L}^L \ge \sum_{\forall m} A^L_{m,c} \mathcal{T}_{m,c}, \forall c \in C^L
\label{working_hrs_line}
\end{equation}
\begin{equation}
\mathcal{L}^T \ge \sum_{\forall m} A^T_{m,c} \mathcal{T}_{m,c}, \forall c \in C^T
\label{working_hrs_tree}
\end{equation}
}
\begin{equation}
\sum_{\forall c \in C^L} A^L_{m,c} = 1, \forall m \in \Omega_{DK}
\label{Pre_assignL}
\end{equation}
\begin{equation}
\sum_{\forall c \in C^T} A^T_{m,c} = 1, \forall m \in \Omega_{DK}
\label{Pre_assignT}
\end{equation}
\begin{equation}
\sum_{\forall r} Cap^R_r Res^C_{c,w,r} \le (\delta_{w,c}+z_{w,c}) Cap^C_c, \forall w,c \in C^L
\label{Pre-Capacity}
\end{equation}
\vspace{-0.2cm}
\begin{equation}
 z_{w,c} \le \delta_{w,c} , \forall w,m,c \in C^L
\label{z<delta}
\end{equation}
\begin{equation}
\mathcal{P}_{c,w} \ge A^L_{m,c} tr_{w,m} - M (1-z_{w,c}), \forall w,m,c \in C^L
\label{Penalty_distance}
\end{equation}
\begin{equation}
\sum_{\forall c \in C^L} Res^C_{c,w,r} \le Res^D_{w,r}, \forall w,r
\label{ResC<ResD}
\end{equation}
\begin{equation}
\sum_{\forall w} Res^C_{c,w,r} \ge \sum_{\mathclap{\forall m}} A^L_{m,c} \mathcal{R}_{m,r}, \forall c \in C^L,r
\label{ResC>AL_R}
\end{equation}
{\color{black}
\begin{equation}
\bar{tr} \ge tr_{m,n} (A^L_{m,c}+A^L_{n,c}-1), \forall m,n,c \in C^L
\label{PreDist1}
\end{equation}
\begin{equation}
\bar{tr} \ge tr_{w,m} (\delta_{w,c}+A^L_{m,c}-1), \forall w,m,c \in C^L
\label{PreDist2}
\end{equation}
\begin{equation}
\bar{tr} \ge tr_{m,n} (A^T_{m,c}+A^T_{n,c}-1), \forall m,n,c \in C^T
\label{PreDist3}
\end{equation}
\begin{equation}
\bar{tr} \ge tr_{w,m} (\delta_{w,c}+A^T_{m,c}-1), \forall w,m,c \in C^T
\label{PreDist4}
\end{equation}
}
\begin{equation}
\{A^{L/T},z\} \in \{0,1\}, \{\mathcal{P}, Res^C\} \ge 0
\label{Assign_var_def}
\end{equation}
}
{\color{black}The objective (\ref{Pre_obj}) consists of four parts. The first two terms minimize the expected time of the last repair for the line crews ($\mathcal{L}^L$) and tree crews ($\mathcal{L}^T$). The variables $\mathcal{L}^L$ and $\mathcal{L}^T$ are defined in constraints (\ref{working_hrs_line}) and (\ref{working_hrs_tree}), respectively. The third term in (\ref{Pre_obj}) is a penalty cost used to limit the number of times a crew goes back to the depot to pick up additional resources. The fourth term $\bar{tr}$ is the maximum travel time for the crews.} Constraints (\ref{Pre_assignL})-(\ref{Pre_assignT}) assign each damaged component to one crew. The amount of resources a crew can carry is limited by the crew's capacity in (\ref{Pre-Capacity}). Binary variable $z_{w,c}$ is equal to 1 if a crew requires additional resources. In such case, the crew goes back to the depot to pick up the required resources. Constraint (\ref{z<delta}) states that the crews can go back to the depot they started from. We set the penalty term $\mathcal{P}_{w,c}$ to be equal to the maximum travel time between the depot and the assigned damage components, as defined in (\ref{Penalty_distance}). The big $M$ constant is added so that the penalty term equals 0 if the crew does not go back to the depot for additional resources. The crews must use the resources available in the depot as enforced by (\ref{ResC<ResD}). Constraint (\ref{ResC>AL_R}) indicates that the number of resources crew $c$ has should be enough to repair the assigned damaged components. {\color{black}Constraints (\ref{PreDist1})-(\ref{PreDist4}) are used to identify the maximum travel time between the damaged components that are assigned to each crew. If components $m$ and $n$ are assigned to crew $c$, then $\bar{tr} \ge tr_{m,n}$.}

\subsubsection{Initial Solution and Optimization}

After assigning each damaged component to a crew, DSRRP is solved with the crews dispatched to the assigned components. Subsequently, a neighborhood search method is used to improve the initial route. The optimization problem considered in this paper involves a dynamically changing environment due to the uncertainty of the repair time, solar irradiance, and demand. The repair time is updated periodically either by the repair crews or the damage assessors. Therefore, we apply the neighborhood search algorithm continuously and update the routing solution as more information is obtained. The advantage of this method is that it allows the algorithm to update the solution while the repair crews are repairing the lines, therefore, loosening the time limit restriction. The pseudo-code for the proposed algorithm, referred to as the Reoptimization algorithm, is detailed in Algorithm 1. 

\begin{algorithm}[h!]
 \caption{Reoptimization Algorithm for DSRRP}
 \small
 \begin{algorithmic}[1]
 \renewcommand{\algorithmicrequire}{\textbf{Input:}}
 \renewcommand{\algorithmicensure}{\textbf{Output:}}
 
  \STATEx Obtain the location of the outages from the damage assessors.
  \STATE solve using \textbf{CPLEX} \COMMENT{Assignment}
  \STATEx $(A^L,A^T) = \argmin \{\eqref{Pre_obj}|{\rm{s.t.}~}\eqref{working_hrs_line}\mbox{-}\eqref{Assign_var_def}\}$
  \FORALL {$c \in C^L$}
  \STATE $N(c) = \{m|\forall m \in \Omega_{DK}, A^L_{m,c}=1\} \cup \Omega_P$
  \ENDFOR
  \FORALL {$c \in C^T$}
  \STATE $N(c) = \{m|\forall m \in \Omega_{DT}, A^T_{m,c}=1\} \cup \Omega_P$
  \ENDFOR
  \STATE solve using \textbf{CPLEX} (time limit = 300 s) \COMMENT{Assignment-DSRRP} 
  \STATEx $\zeta^* = \min\{\eqref{DSRRP_Obj}|{\rm{s.t.}~}\eqref{P_clpu}\mbox{-}\eqref{DSRRP_var_def},\sum_{n \in {N}(c)}x_{m,n,c}=1, \forall c, m \in {N}(c)\}$
 \STATE obtain solution $x^*$ and objective $\zeta^*$
 \STATE let $\bar{x}=x^*$ and $\bar{\zeta}=\zeta^*$
 \REPEAT
 \STATE set $count = 0$
\STATE set $ss = ss_0$ \COMMENT{sample size}
 \WHILE {time limit is not surpassed}
 \COMMENT{Neighborhood Search} 
 \STATE let $\bar{N} = sample(N,ss)$, where $\bar{N} \subset N$ and $|\bar{N}|=ss$.
 \STATE solve using \textbf{CPLEX} {\color{black}(time limit = 120 s)} with \textit{warm start} 
 \STATEx \hspace{\algorithmicindent}\hspace{\algorithmicindent}{$\zeta^* =~ \min\{\eqref{DSRRP_Obj}|{\rm{s.t.}~}\eqref{P_clpu}\mbox{-}\eqref{DSRRP_var_def},x_{m,n,c}~=~\bar{x}_{m,n,c}, \forall c,~ m \in$
  \STATEx \hspace{\algorithmicindent}\hspace{\algorithmicindent}${N}\backslash\bar{{N}}, n \in {N}\backslash\bar{{N}}\}$}
 \STATE obtain $x^*$ and objective $\zeta^*$
 \IF {$\zeta^* < \bar{\zeta}$}
 \STATE set $\bar{x}=x^*; \bar{\zeta}=\zeta^*; count = 0$
 \ELSE
 \STATE $count = count + 1$
 \ENDIF
 \STATE \textbf{if} {$ss = |N|$} \textbf{then} \textbf{break} \COMMENT{solution is optimal}
 \STATE \textbf{if} {$count = h_1$} \textbf{then} $ss = ss + 1$ 
 \STATE \textbf{if} {$count = h_1+h_2$} \textbf{then} \textbf{break} 
 \ENDWHILE
\STATE dispatch crews and set the traveled path as constant
\STATE {update the repair time and return to Step 11}
\UNTIL{all lines are repaired}
 \end{algorithmic} 
 \end{algorithm}

In Step 1, the assignment problem is solved using CPLEX \cite{CPLEX_2016} to obtain the binary variables $A^L_{m,c}$ and $A^T_{m,c}$. These variables are used to find $N(c)$, which is the set of damaged components assigned to crew $c$. For example, consider the set of damaged components $\Omega_{DK}= \{1,2,3,4,5\}$, if line crew 1 is assigned with damaged components 1 and 3, then $A^L_{m,c} = \{1,0,1,0,0\}$ and $N(1) = \{1,3\} \cup \Omega_{P}$. $N(c)$ is found for each crew in Steps 2-7. Consequently, a simplified DSRRP is solved in Step 8 by allowing the crews to only repair the assigned damaged components. In Step 10, the obtained route $x^*$ and objective $\zeta^*$ are set to be the incumbent (current best solutions) route ($\bar{x}$) and objective ($\bar{\zeta}$). 

\begin{figure}[h!]
\setlength{\abovecaptionskip}{0pt} 
\setlength{\belowcaptionskip}{0pt} 
  \centering
    \includegraphics[width=0.45\textwidth]{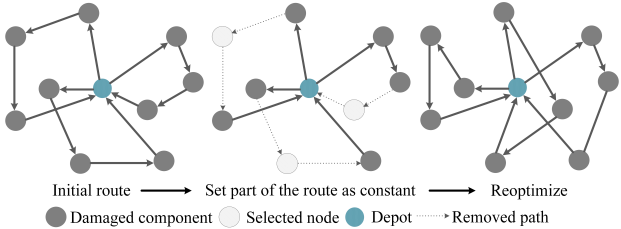}
  \caption{A single iteration of the neighborhood search, with $|\bar{N}|=3$.}\label{Reopt_alg}
\end{figure}

Steps 11-29 represent the neighborhood search algorithm. The algorithm selects a subset of damaged components $\bar{N}$, where $\bar{N} \subset N$, then removes the paths connected to $\bar{N}$ and sets the rest of the routes to be constant by forcing $x_{m,n,c}~=~\bar{x}_{m,n,c}, \forall c,~ m \in {N}\backslash\bar{{N}}, n \in {N}\backslash\bar{{N}}$. Afterwards, DSRRP is solved to obtain an improved solution, the process is demonstrated in Fig. \ref{Reopt_alg}, where $|\bar{N}|=3$. 
 
Steps 12 and 13 initialize a counter and the sample size ($ss$), respectively. In Step 15, the subset $\bar{N}$ is determined by randomly selecting $ss$ nodes from $N$. The parameters $ss_0$, $h_1$, and $h_2$ are constants used to tune the algorithm. The value of $ss_0$ determines the size of the subset $\bar{N}$ in the first iteration. The size of $\bar{N}$ is increased after $h_1$ iterations with no change to the objective, and the neighborhood search algorithm is terminated after $h_1+h_2$ iterations with no change to the objective. In this paper, $ss_0$ is set to be 3, as selecting 1 damaged component will not change the route, and selecting 2 has minimal impact on the route. The values of $h_1$ and $h_2$ were determined experimentally using several test cases, both $h_1$ and $h_2$ equal 3. 

\begin{figure}[h!]
\setlength{\abovecaptionskip}{0pt} 
\setlength{\belowcaptionskip}{0pt} 
  \centering
    \includegraphics[width=0.45\textwidth]{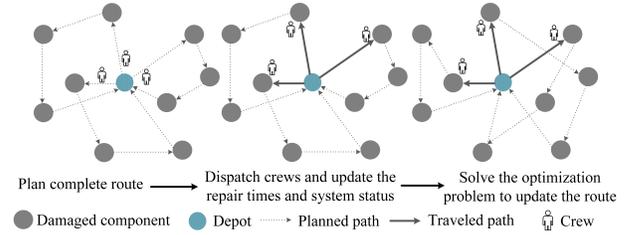}
  \caption{Dynamic vehicle routing problem.}\label{ResLearn}
\end{figure}

\begin{figure}[h!]
\setlength{\abovecaptionskip}{0pt} 
\setlength{\belowcaptionskip}{0pt} 
  \centering
    \includegraphics[width=0.4\textwidth]{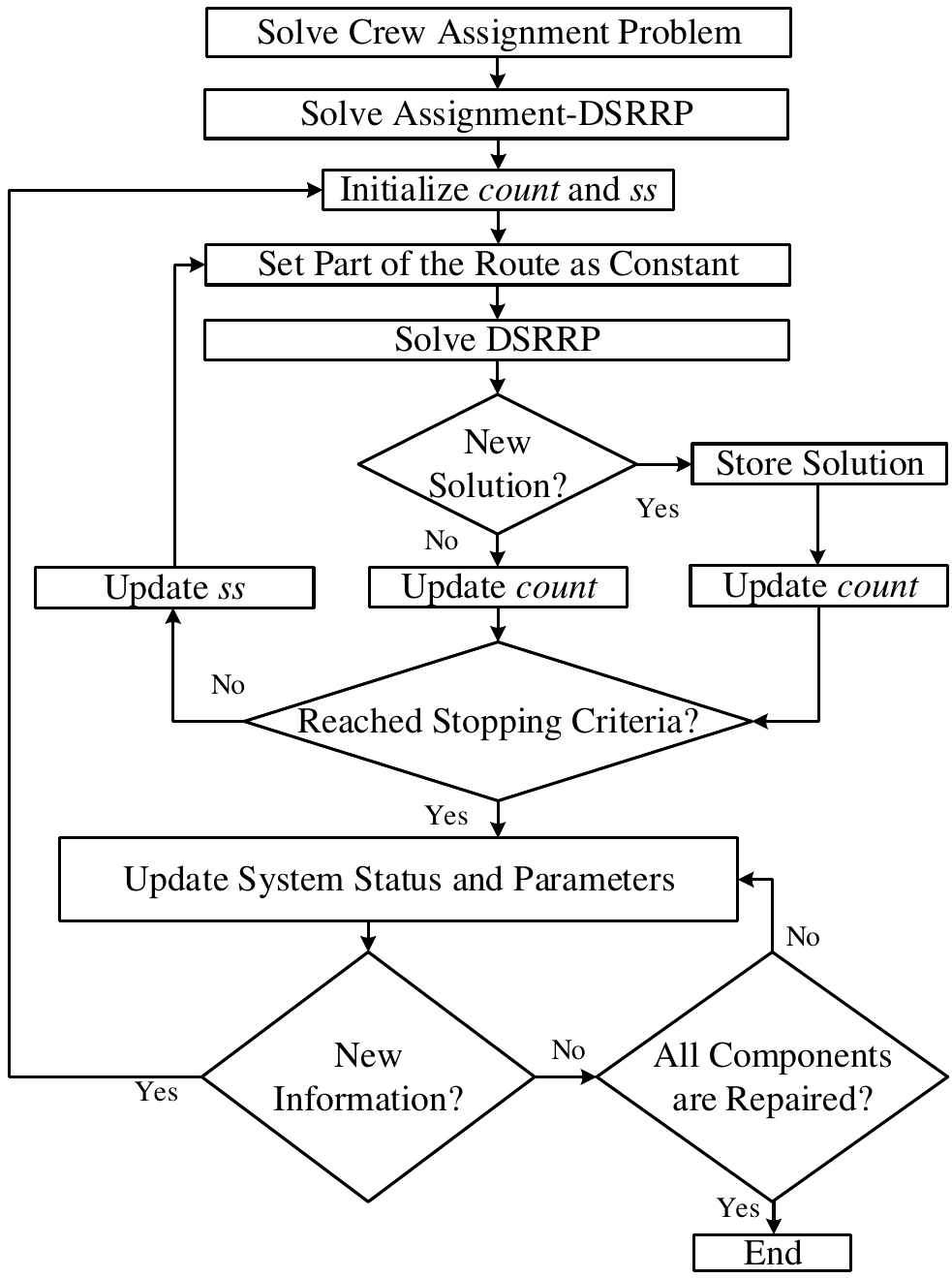}
  \caption{{\color{black}Flow chart of the Reoptimization algorithm.}}\label{reopt_flowchart}
\end{figure}

The DSRRP is solved in Step 16 with parts of the route set as constant. {\color{black}To obtain a fast solution, we warm start (provide a starting point) CPLEX by using the incumbent solution and enforce a time limit of 120 seconds for each iteration.} The objective value $\zeta^*$ obtained from Step 16 is compared to the current incumbent solution $\bar{\zeta}$. If the value is improved, we set $\zeta^*$ and $x^*$ as the current incumbent solutions and update the counter, otherwise, the counter increases by one. The process is repeated until the counter reaches $h_1$, where we increase the size of the subset in Step 24. If the sample size is $|N|$; i.e., the complete problem is solved without simplification, then the solution is optimal and the neighborhood search stops. Also, the search ends once the counter reaches $h_1+h_2$, or if the time limit is reached. The crews are then dispatched to the damaged components, and the traveled paths are set as constants in the optimization problem. After that, the repair time is updated and Steps 14-26 are repeated to update the route, as shown in Fig. \ref{ResLearn}. The idea of the dynamic approach is to run Steps 14-26 while maintaining the best solution in an adaptive memory. Once the operator receives an update from the field, the neighborhood search is restarted with the newly acquired information. Whenever a crew finishes repairing the assigned damaged component, the crew is provided with the current best route $\bar{x}$. A flowchart for the proposed algorithm is presented in Fig. \ref{reopt_flowchart}.

\subsection{Priority-based Method} \label{PB_method}
In general, utilities schedule the repair using a defined
restoration priority lists. To compare the proposed approach to current practices, a priority-based method is developed to replicate the procedure that the utilities follow. Each utility has its own priority list but it can be generally summarized as follows \cite{Utilities_Priorities}:
\begin{enumerate}
\item Repair lines connected to high-priority customers. 
\item Repair three-phase lines starting with upstream lines
\item Repair single phase lines and individual customers  
\end{enumerate}
Define $L_p$ as the set of lines to repair with priority $p$, and $w_p$ is a weighting factor, where $w_1>w_2>w_3$ (e.g., $w_1=10, w_2=5, w_3=1$). $L_1$ contains the lines that must be repaired to restore critical customers, $L_2$ represents the three-phase lines not in $L_1$, and $L_3$ represents the rest of the lines. The following routing model is solved to find the repair schedule by utilizing the priority of each line, as follows:
{\small
\begin{equation}
x^p \mbox{=} \argmin \{\sum_{\forall p}\sum_{\forall k \in L_p}\sum_{\forall c \in C^L}w_p \alpha_{c,k} |{\rm{s.t.}} (23)\mbox{-}(38)\}
\label{priority_model}
\end{equation}
}
The objective of (\ref{priority_model}) is to minimize the arrival time of the line crews at each damaged components, while prioritizing the high-priority lines through multiplying the arrival time by the weight $w_p$. The priority-based model is similar to DSRRP, but without the power operation constraints. However, it is still difficult to solve directly in a short time using a commercial solver such as CPLEX. Therefore, the same procedure presented in Algorithm 1 is used to solve (\ref{priority_model}). After obtaining the route $x^p$, the DSRRP problem is solved by setting $x=x^p$; i.e., we solve $\min\{(1)|~{\rm{s.t.}~}(2)\mbox{-}(40),x_{m,n,c}~=~x^p_{m,n,c}, \forall c, m, n$\}.
\section{Simulation and Results}\label{chap:6}

Modified IEEE 123- and 8500-bus distribution feeders are used as test cases for the DSRRP problem. Detailed information on the networks can be found in \cite{123ieee} and \cite{8500ieee}. {\color{black}Since transportation networks data for the IEEE 123- and 8500-bus test cases are not available, the network $G$ and the travel times are simulated by using the Euclidean distance \cite{Hent2015}. The average speed of the crews is assumed to be 35 mph in the simulated problems. The travel time is calculated by dividing the Euclidean distances between all nodes by the speed of the crews. We then scale the travel time such that the travel time between the two furthest locations equals 2 hours. The x-coordinates and y-coordinates for the IEEE 123- and 8500-bus test cases can be found in \cite{123ieee} and \cite{8500ieee}, respectively. We assume that there is an available path to each damaged component.}.

\begin{figure}[h!]
\setlength{\abovecaptionskip}{0pt} 
\setlength{\belowcaptionskip}{0pt} 
  \centering
    \includegraphics[width=0.45\textwidth]{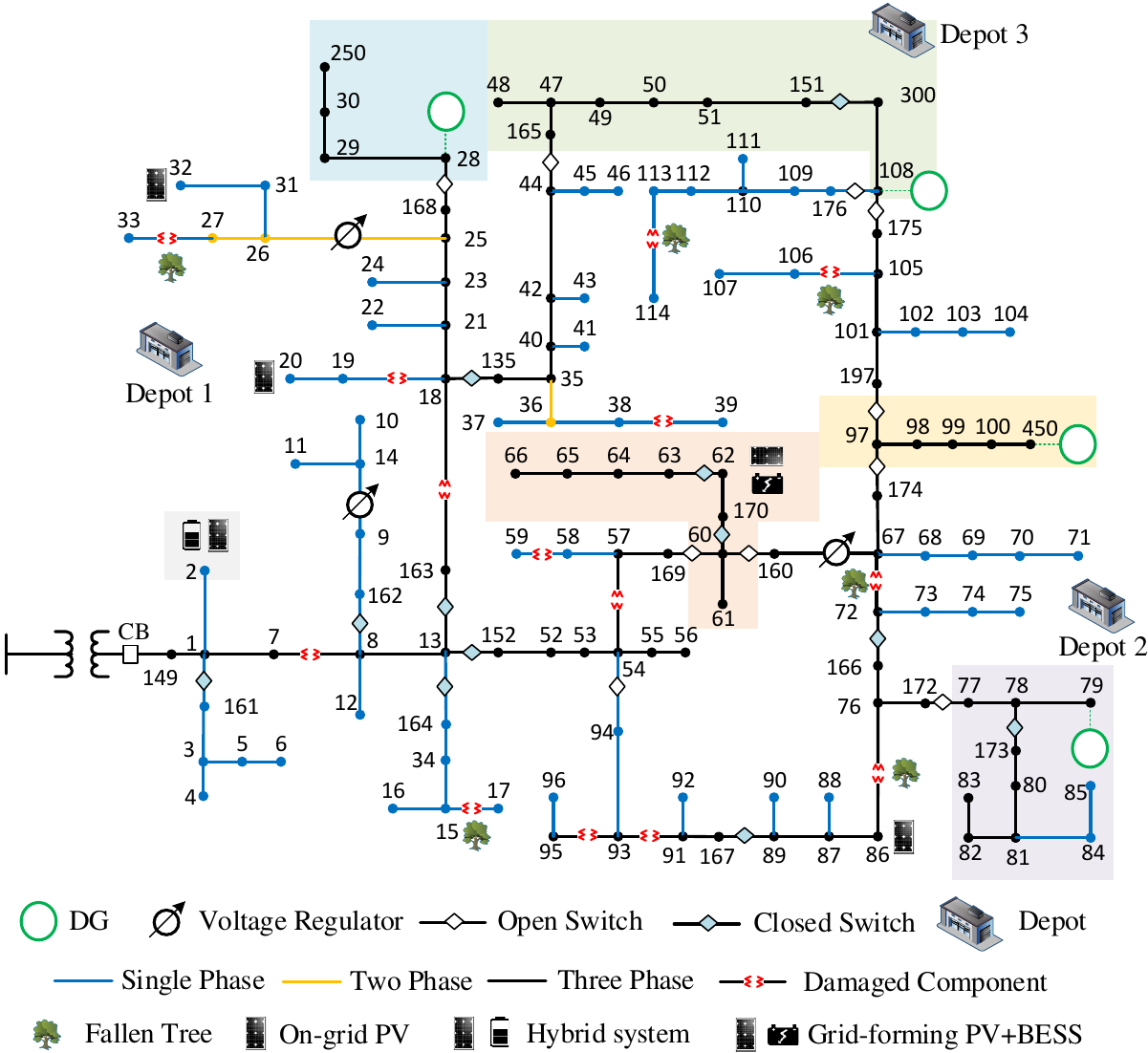}
  \caption{{\color{black}Initial state of the distribution network after 14 lines are damaged.}}\label{IEEE123_0}
\end{figure}

\begin{figure}[h!]
\setlength{\abovecaptionskip}{0pt} 
\setlength{\belowcaptionskip}{0pt} 
  \centering
    \includegraphics[width=0.4\textwidth]{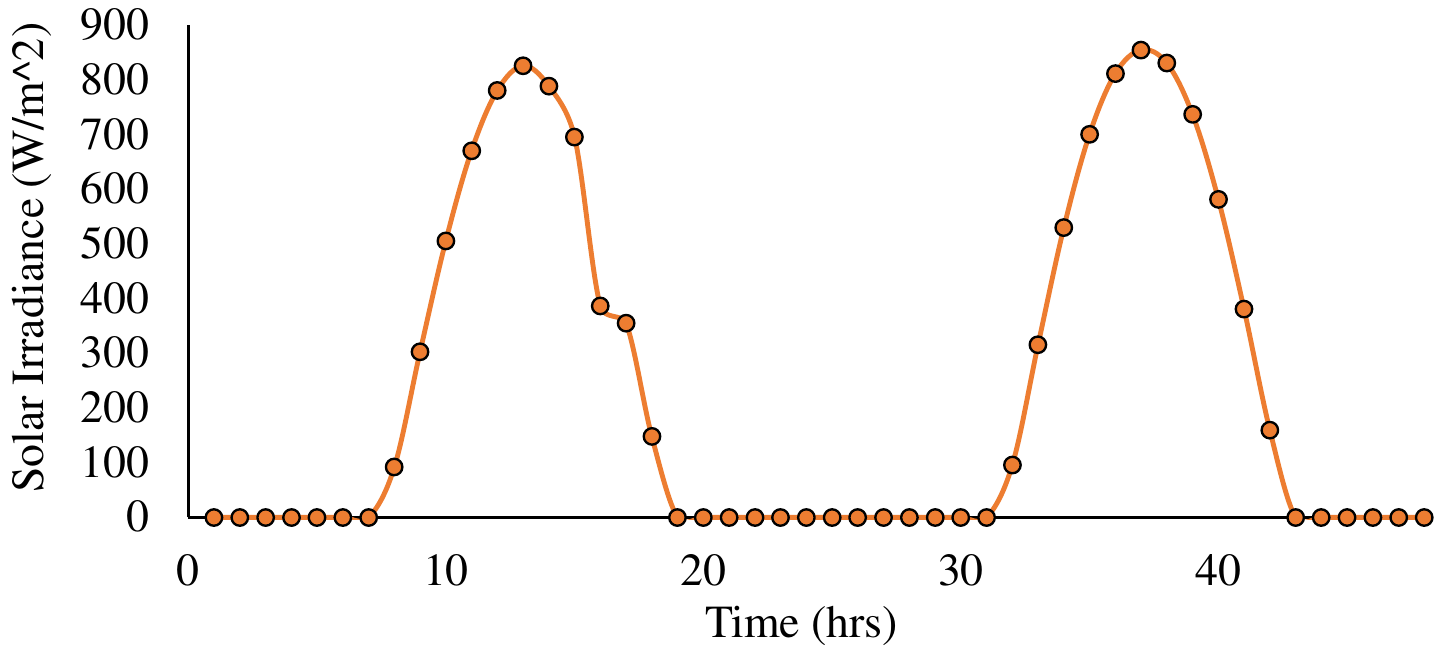}
  \caption{{\color{black}Solar irradiance for the PV systems in the simulation \cite{Solar_Data}.}}\label{SolarShape}
\end{figure}

\begin{figure}[h!]
\setlength{\abovecaptionskip}{0pt} 
\setlength{\belowcaptionskip}{0pt} 
  \centering
    \includegraphics[width=0.49\textwidth]{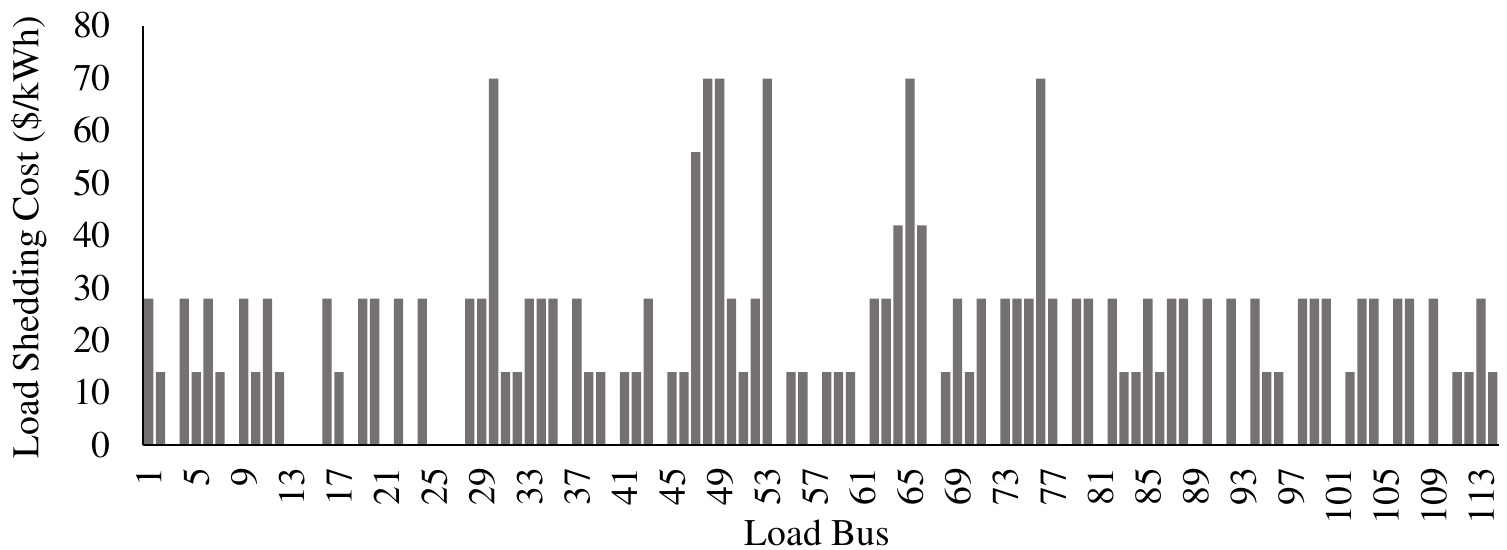}
  \caption{{\color{black}The load shedding cost in \$/kWh of each load in the simulation.}}\label{LoadCostFig}
\end{figure}

The IEEE 123-bus feeder, shown in Fig. \ref{IEEE123_0}, is modified by including 4 dispatchable DGs, 18 new switches, 5 PVs and 2 BESSs. The 4 DGs are rated at 300 kW and 250 kVAr. PVs in On-grid and hybrid systems are rated at 50 kW, and the PV at bus 62 is rated at 900 kW. The forecasted solar irradiance used in the simulation is presented in Fig. \ref{SolarShape}, which is obtained from the National Solar Radiation Data Base (NSRDB) \cite{Solar_Data}. The data in Fig. \ref{SolarShape} represent the solar irradiance at a location impacted by Hurricane Matthew. The BESSs at bus 2 and 62 are rated at 50 kW/132 kWh and 500 kW/ 2100 kWh, respectively. 
Fig. \ref{LoadCostFig} shows the load shedding costs of each load. The problems of optimally allocating the resources, DGs, or switches, are out of the scope of this paper. We assume there are 3 depots, 6 line crews distributed equally between the depots, and 4 tree crews with 2 located in Depot 2 and and 1 tree crew in each of the other depots. The time step in the simulation is 1 hour. The simulated problem is modeled in AMPL and solved using CPLEX 12.6.0.0 on a PC with Intel Core i7-4790 3.6 GHz CPU and 16 GB RAM.


\vspace{-0.3cm}
\subsection{DSRRP solution comparison}

The repair and restoration problem is solved using five methods: 1) a cluster-first DSRRP-second (C-DSRRP) approach presented in \cite{Arif2016}, the method clusters the damaged components to the depot, then solves DSRRP; 2) the priority-based method presented in Section \ref{PB_method}; 3) an assignment-based method where the damaged components are assigned to the crews, then DSRRP is solved (A-DSRRP), which is similar to Steps 1-8 in Algorithm 1; 4) Reoptimization algorithm; 5) CPLEX with warm start using the Reoptimization algorithm solution. 

Once an outage occurs, the distribution network is reconfigured, and the DGs are dispatched to restore as many customers as possible, before conducting the repairs. A random event is generated on the IEEE 123-bus system, where 14 lines are damaged, four of which are damaged by trees. Fig. \ref{IEEE123_0} shows the recovery operation of the distribution system to the outages before the repairs; i.e., the state of the system at time $t=0$. The solution shown in Fig. \ref{IEEE123_0} is obtained regardless of the solution algorithm used, as the algorithms will only affect the repair schedule and the network operation during the repairs. Before the outage, all switches are closed except 151-300 and 54-94. {\color{black}Since line 7-8 is damaged, the circuit breaker at the substation is opened. Sectionalizer 28-168 is switched off, forming a small microgrid, to serve the loads at buses 28 to 30. Similarly, switches 44-165, 77-172, 97-174, 97-197, 108-175 and 108-176 are opened and 151-300 is closed to form additional microgrids using the DGs in the network. Switches 60-160 and 60-169 are opened so that the PV+BESS at bus 62 can form a microgrid. The battery at bus 2 can serve the local load in the first few hours after the damage.} The repair/tree-clearing times and required resources are given in Table \ref{Data_input}. The estimated repair time is assumed to be accurate. It is assumed that each crew can carry 30 units of resources, and the required capacities ($Cap^R_r$) for the 6 types of resources are \{3, 2.5, 2, 1, 4, 1\}. A summary of the results and performances of different solution methods is shown in Table \ref{summary_alg}. The time limit is set to be 3600 seconds \cite{last_mile} for all methods except for the last one (CPLEX with a warm start) in order to find the optimal solution.

\begin{table}[htbp]
  \centering
  \footnotesize
  \caption{The resources and time required to repair the damages}
  \vspace{-0.2cm}
    \begin{tabular}{|c|c|c|c|c|c|c|c|c|}
    \hline
    \rule{0pt}{2.2ex}\multirow{2}[4]{*}{Line} & \multicolumn{6}{c|}{Resources (units)}                & \multicolumn{2}{c|}{Estimated repair/clearing time (hrs)} \\
\cline{2-9}\rule{0pt}{2.2ex}          & A     & B     & C     & D     & E     & F     & ~~Line Crew~~ & Tree Crew \\
    \hline
    \rule{0pt}{2.2ex}7-8 & 1     & 2     & 0     & 1     & 0     & 0     & 2.5   &  {\cellcolor[gray]{.8}}\\
    \hline
    \rule{0pt}{2.2ex}15-17 & 1     & 2     & 1     & 1     & 0     & 0     & 1.25  & 1 \\
    \hline
    \rule{0pt}{2.2ex}18-19 & 1     & 2     & 1     & 1     & 0     & 0     & 0.5   &  {\cellcolor[gray]{.8}}\\
    \hline
    \rule{0pt}{2.2ex}27-33 & 1     & 2     & 1     & 1     & 0     & 0     & 2.25  &  {\cellcolor[gray]{.8}}\\
    \hline
    \rule{0pt}{2.2ex}38-39 & 1     & 2     & 1     & 1     & 0     & 0     & 1     & 0.75 \\
    \hline
    \rule{0pt}{2.2ex}54-57 & 0     & 2     & 0     & 1     & 2     & 0     & 0.75  &  {\cellcolor[gray]{.8}}\\
    \hline
    \rule{0pt}{2.2ex}58-59 & 1     & 2     & 1     & 1     & 0     & 0     & 0.5   &  {\cellcolor[gray]{.8}}\\
    \hline
    \rule{0pt}{2.2ex}18-163 & 0     & 2     & 0     & 1     & 0     & 2     & 1.75  &  {\cellcolor[gray]{.8}}\\
    \hline
    \rule{0pt}{2.2ex}67-72 & 0     & 2     & 0     & 1     & 0     & 0     & 4     & 1.25 \\
    \hline
    \rule{0pt}{2.2ex}76-86 & 1     & 2     & 1     & 1     & 0     & 0     & 6     & 2 \\
    \hline
    \rule{0pt}{2.2ex}91-93 & 0     & 2     & 0     & 1     & 2     & 0     & 1.5   &  {\cellcolor[gray]{.8}}\\
    \hline
    \rule{0pt}{2.2ex}93-95 & 1     & 2     & 1     & 1     & 0     & 0     & 2.75  &  {\cellcolor[gray]{.8}}\\
    \hline
    \rule{0pt}{2.2ex}105-106 & 1     & 2     & 1     & 1     & 0     & 0     & 1.75  & 1 \\
    \hline
    \rule{0pt}{2.2ex}113-114 & 1     & 2     & 1     & 1     & 0     & 0     & 0.75  & 0.5 \\
    \hline
    \end{tabular}%
  \label{Data_input}%
\end{table}%

\begin{table}[htbp!]
{\color{black}
\footnotesize
  \centering
  \caption{A comparison between four methods for the IEEE 123-bus system}
  \vspace{-0.2cm}
    \begin{tabular}{cccccc}
    \hline
   \rule{0pt}{2ex}\multirow{2}[2]{*}{Method} & Objective & Optimality & CPU & Load & Restoration \\
          &  Value & Gap & Time  &   Served    & Time  \\
    \hline
    \rule{0pt}{2ex} C-DSRRP &    \$241,371   &   21.16\% & 3600 s    &    61.86 MWh&  12 hrs\\
    Priority-based &    \$229,112   &   15.01\% & 662 s    &    62.25 MWh &  9 hrs\\
    A-DSRRP & \$211,597 & 6.21\% & 206 s & 62.98 MWh & 9 hrs \\
    Reoptimization & \$199,210 & 0.00\% & 694 s & 63.5 MWh& 9 hrs \\
    CPLEX& \$199,210 & 0.00\% & 4 hrs & 63.5 MWh& 9 hrs \\
    \hline
    \end{tabular}%
  \label{summary_alg}%
}
\end{table}%

\begin{figure}[h!]
\setlength{\abovecaptionskip}{0pt} 
\setlength{\belowcaptionskip}{0pt} 
  \centering
    \includegraphics[width=0.49\textwidth]{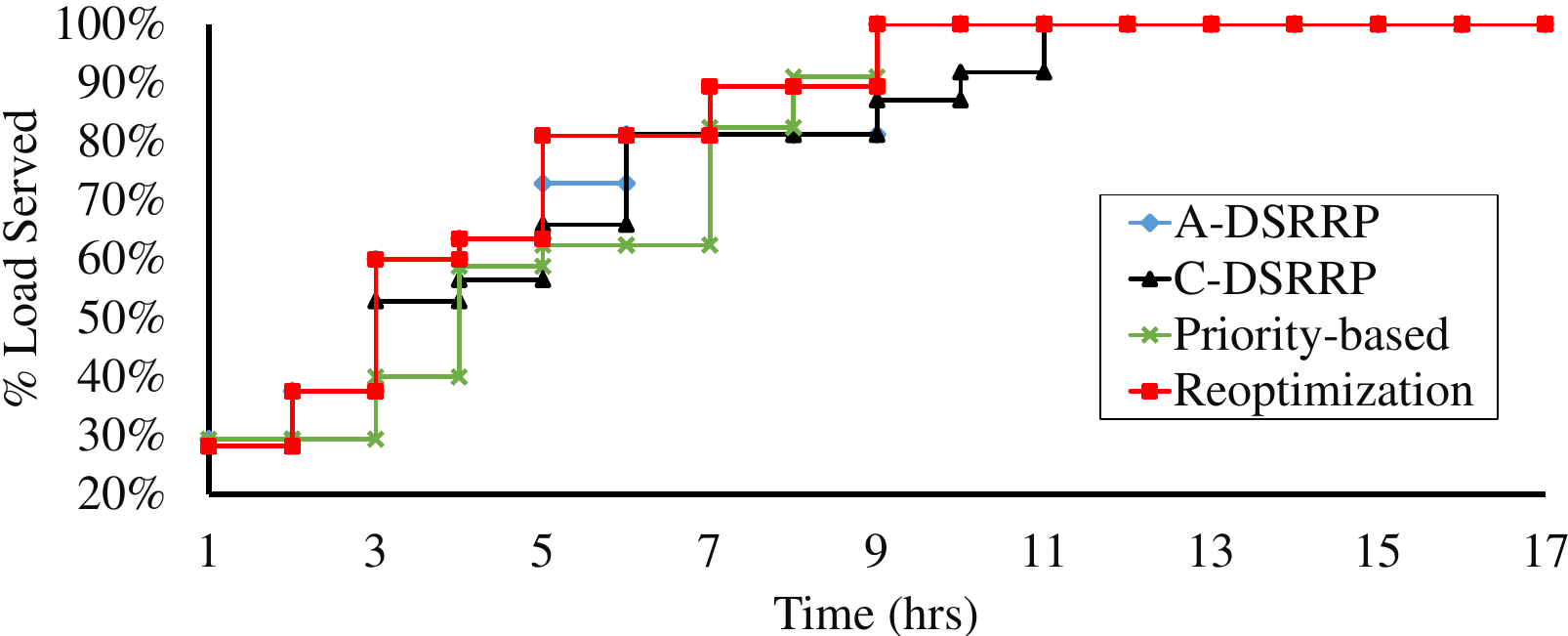}
  \caption{Percentage of load served at each time step.}\label{Load_Served_fig}
\end{figure}

The fifth column in Table \ref{summary_alg} is the amount of energy served, and the sixth column (restoration time) is the time when all loads are restored. The assignment-based approach (A-DSRRP) is the fastest but the solution is not optimal, neighborhood search in the Reoptimization algorithm improved the routing solution and obtained the best repair schedule. {\color{black}To obtain the optimal solution, the route obtained from the proposed method is used to warm start CPLEX and solve DSRRP. CPLEX showed that the solution obtained from the Reoptimization algorithm is optimal. C-DSRRP reached the time limit but produced a feasible solution with 21.16\% optimality gap, while the priority-based method achieved an objective value which is \$29,902 higher than the optimal solution.} The change in percentage of load served for each method is shown in Fig. \ref{Load_Served_fig}. The proposed algorithm outperformed the other methods.

{\color{black}
Next, we compare the Reoptimization algorithm with the priority-based method using three different damage scenarios on the IEEE 123-bus system. The simulation results are shown in Table \ref{3test_123}. The proposed method outperforms the priority-based method in all instances with comparable computation times. The results show how the proposed algorithm can achieve near-optimal solutions, and indicate the importance of co-optimizing repair scheduling and the operation of the distribution system. For the first test case, the algorithm achieved the optimal solution, while the optimality gap for the priority-based method is 2.98\%. The Reoptimization algorithm achieved solutions that are approximately 11\% and 17\% less than the priority-based method for the second and third test cases, respectively.

\begin{table}[htbp]
  \centering
  \footnotesize
  \caption{Three test cases solved using the Reoptimization and priority-based methods}
    \vspace{-0.2cm}
    \begin{tabular}{|c|c|c|c|c|c|c|}
    \hline
    \multirow{2}[4]{*}{Damage} & \multicolumn{3}{c|}{Reoptimization} & \multicolumn{3}{c|}{Priority-based} \\
\cline{2-7}                & Obj.  & \% Gap & CPU Time & Obj.  & \% Gap & CPU Time \\
    \hline
    15 Lines & \$158,023 & 0.00\% & 660 s & \$162,734 & 2.98\% & 464 s \\
    \hline
    20 Lines & \$248,986 & 2.53\% & 762 s & \$279,197 & 14.97\% & 392 s \\
    \hline
    25 Lines & \$388,760 & 2.27\% & 782 s & \$467,278 & 22.93\% & 520 s \\
    \hline
    \end{tabular}%
  \label{3test_123}%
\end{table}%
}
\subsection{Dynamic DSRRP}

In practice, the crew repair time is continuously changing. Moreover, the dispatch commands must be issued as fast as possible to reduce the outage duration. Therefore, the DSRRP must be solved efficiently and the solutions should be dynamically updated according to the current crew repair time. 
To simulate the change in repair time, it is assumed that once a crew reaches the damaged component, the repair time is updated to its actual value by adding a random number from the continuous uniform distribution on [-2,2] to the estimated time. For example, once crew 1 arrives at line 7-8, the repair time is changed from 2.5 to 3 hours. {\color{black}Similarly, the solar irradiance is updated by adding $\pm$5\% to the forecasted value}. The time limit at Step 14 in Algorithm 1 is set to be 15 minutes after the first dispatch, so that the repair time is updated every 15 minutes. While the crews are repairing the damaged components, the neighborhood search algorithm keeps searching for a better solution, and the crews are dispatched using the incumbent solution. 

\begin{figure}[h!]
\setlength{\abovecaptionskip}{0pt} 
\setlength{\belowcaptionskip}{0pt} 
  \centering
    \includegraphics[width=0.49\textwidth]{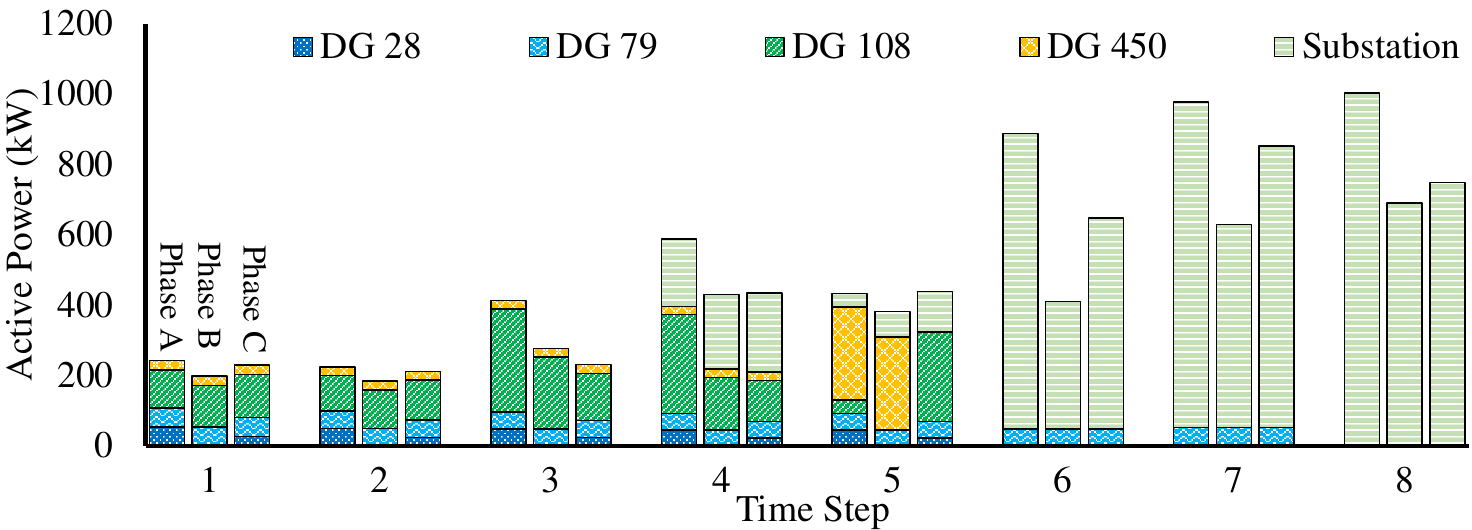}
  \caption{The 3-phase active power delivered by the DGs and substation.}\label{DG_output}
\end{figure}

\begin{figure}[h!]
\setlength{\abovecaptionskip}{0pt} 
\setlength{\belowcaptionskip}{0pt} 
  \centering
    \includegraphics[width=0.49\textwidth]{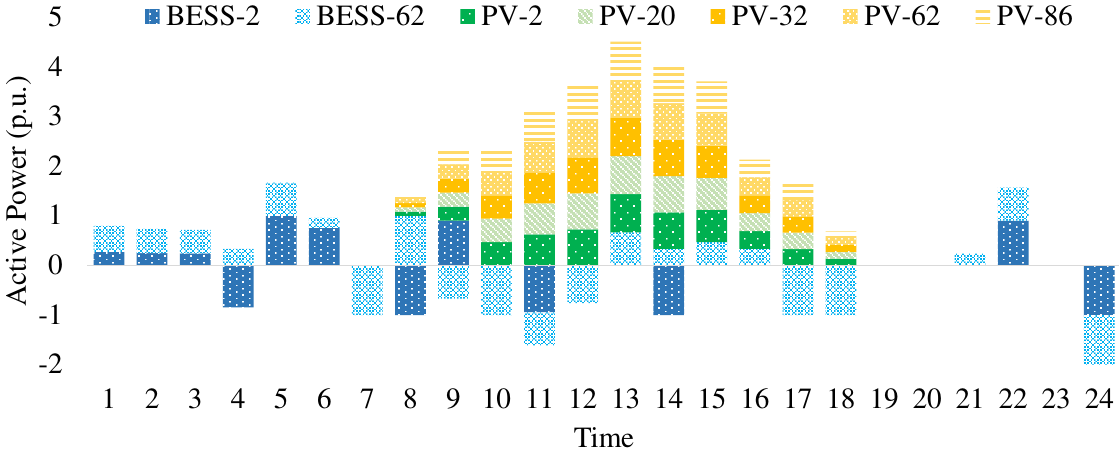}
  \caption{{\color{black}The active power delivered by the PVs and BESSs.}}\label{PV_BESS_Fig}
\end{figure}

The complete route is given in Table \ref{RH_route_table}. The total cost is \$192,694, and the total energy served is 64.7 MWh. Table \ref{Event_timeline} shows the timeline of events after solving DSRRP, where all loads are restored after 8 hours. The initial states of the switches are shown in Fig. \ref{IEEE123_0}, and the subsequent switching operations are given in Table \ref{Event_timeline}. The 3-phase output of the DGs and the substation are shown in Fig. \ref{DG_output}, and {\color{black}Fig. \ref{PV_BESS_Fig} shows the output of the PVs and BESSs. Crew 5 repairs line 38-39 and switch 18-135 is opened and 44-165 is closed to restore the loads at buses 35 to 46. Once line 113-114, is repaired by tree crew 10 and line crew 4, switch 108-174 is closed to restore the loads at buses 109 to 114.  After repairing line 7-8 in time step 4, the CB is closed and the network starts to receive power from the substation. Switches 13-163 and 13-164 are opened to keep lines 15-17, 18-19, and 27-33 isolated. Loads at buses 52 to 59 are restored after repairing lines 54-57 and 58-59. 8 loads are restored after repairing lines 15-17 and 105-106. After 6 hours, the loads around depot 1 are restored after repairing line 18-19 and closing switch 13-163. Finally, all loads are restored after 8 hours once lines 76-86, 91-93, and 93-95 are repaired.} Switch 151-300 is opened and 18-135 is closed to return the network to its original configuration, and the substation can serve all loads.   
\begin{table}[h!]
\footnotesize
  \centering
  \caption{Routing solution for the dynamic 123-bus test case}
  \vspace{-0.2cm}
\renewcommand{\tabcolsep}{0.3mm}
    \begin{tabular}{l|ccccccccccc}
    \hline
     \rule{0pt}{2ex}Crew & \multicolumn{9}{c}{Route}\\
    \hline
 \rule{0pt}{2ex}Crew 1 & DP 1 &$\rightarrow$& 7-8  & $\rightarrow$ & 15-17 &  &  &  &  & &  \\
 Crew 2 & DP 1 & $\rightarrow$& 163-18  & $\rightarrow$ & 27-33 & $\rightarrow$ & DP 1 & $\rightarrow$ & 93-95 & &  \\
 Crew 3 & DP 2 & $\rightarrow$& 54-57 & $\rightarrow$ & 18-19 &  &  &  & & &  \\
 Crew 4 & DP 2 & $\rightarrow$& 113-114 & $\rightarrow$ & DP 3 & $\rightarrow$ & 105-106 & $\rightarrow$ & DP 2 & $\rightarrow$ & 91-93 \\
 Crew 5 & DP 3 & $\rightarrow$& 38-39 & $\rightarrow$ & 67-72 &  & &  &  &  &  \\
 Crew 6 & DP 3 & $\rightarrow$& 58-59 & $\rightarrow$ & 76-86 &  &  &  &  & &  \\
 Crew 7 & DP 1 & $\rightarrow$& 27-33 & $\rightarrow$ & 15-17 &  &  &  &  & &  \\
 Crew 8 & DP 2 & $\rightarrow$& 76-86 &   &  &  &  &  &  & &  \\
 Crew 9 & DP 2 & $\rightarrow$& 67-72 &  &  &  &  &  &  & &  \\
 Crew 10 & DP 3 & $\rightarrow$& 113-114 & $\rightarrow$ & 105-106 &  & &  &  & &  \\
    \hline
    \end{tabular}%
  \label{RH_route_table}%
\end{table}%
\begin{table}[htbp]
\footnotesize
  \centering
  \caption{Event timeline for the IEEE 123-bus dynamic test case}
  \vspace{-0.2cm}
    \begin{tabular}{|c|p{14mm}|p{14mm}|p{18mm}|c|}
    \hline
    \rule{0pt}{2ex}\multirow{1}[4]{*}{Time step} & \multicolumn{2}{c|}{Switch operation} & \multirow{1}[4]{*}{Lines repaired} & \multirow{1}[4]{*}{\% Load Served} \\
\cline{2-3}          & ~~~open & ~~~close &       &  \\
    \hline
    \rule{0pt}{2ex}1     &    {\cellcolor[gray]{.8}}   &   {\cellcolor[gray]{.8}}    &   {\cellcolor[gray]{.8}}    & 29\% \\
    \hline
    \rule{0pt}{2ex}2     &    {\cellcolor[gray]{.8}}    &   {\cellcolor[gray]{.8}}    & {\cellcolor[gray]{.8}} & 29\% \\
    \hline
    \rule{0pt}{2ex}3     & 18-135 & 44-165 \newline 108-176 & 38-39,163-18\newline 58-59,113-114 & 39\% \\
    \hline
    \rule{0pt}{2ex}4     & 13-163\newline 13-164 & 60-169 \newline 150-149 & 7-8\newline 54-57 & 61\% \\
    \hline
    \rule{0pt}{2ex}5     &    {\cellcolor[gray]{.8}}   &   13-164\newline 97-197\newline 108-175    & 15-17\newline 27-33\newline 105-106 & 73\% \\
    \hline
    \rule{0pt}{2ex}6       & 72-166 & 13-163\newline 168-28\newline 60-160 \newline 97-174 & 18-19\newline 67-72 & 89\% \\
    \hline
    \rule{0pt}{2ex}7     &{\cellcolor[gray]{.8}} & {\cellcolor[gray]{.8}} & {\cellcolor[gray]{.8}} & 89\% \\
    \hline
    \rule{0pt}{2ex}8     &    {\cellcolor[gray]{.8}}   & 72-166 \newline 77-172 & 76-86,91-93\newline 93-95 & 100\% \\
    \hline
    \rule{0pt}{2ex}9     & 151-300 & 18-135 &   {\cellcolor[gray]{.8}}    & 100\% \\
    \hline
    \end{tabular}%
  \label{Event_timeline}%
\end{table}%
{\color{black}
\subsection{Algorithm Scalability: IEEE 8500-bus System}
The IEEE 8500-bus feeder test case is used to examine the scalability of the developed algorithm. The test system, shown in Fig. \ref{8500} \cite{8500ieee}, is modified by adding five DGs and five PV systems. The test case has 35 damaged lines, 15 of which are tree induced. We assume there are 3 depots, 12 line crews, and 8 tree crews. The DSRRP problem is solved using the Reoptmization algorithm and the priority-based method. A time limit of 15 minutes is imposed on the algorithms to obtain a solution for dispatching the crews to their first destinations. The total computation time of the priority-based method is 32 minutes (15 for initial dispatch + 17 for updating the routes), and the total computation time for the Reoptimization algorithm is 40 minutes (15 for initial dispatch + 25 for updating the routes). The objective value is 10.2\% lower using the Reoptimization algorithm at \$763,184, compared to \$849,842 when using the priority-based method. Fig. \ref{Load_Served_8500} shows the percentage of load supplied for the two methods. The optimality gap is not known as CPLEX with warm start could not converge to the optimal solution after 24 hours. The simulation results demonstrate the effectiveness of the proposed method and its ability to handle large cases within the time limits.

\begin{figure}[h!]
\setlength{\abovecaptionskip}{0pt} 
\setlength{\belowcaptionskip}{0pt} 
  \centering
    \includegraphics[width=0.4\textwidth]{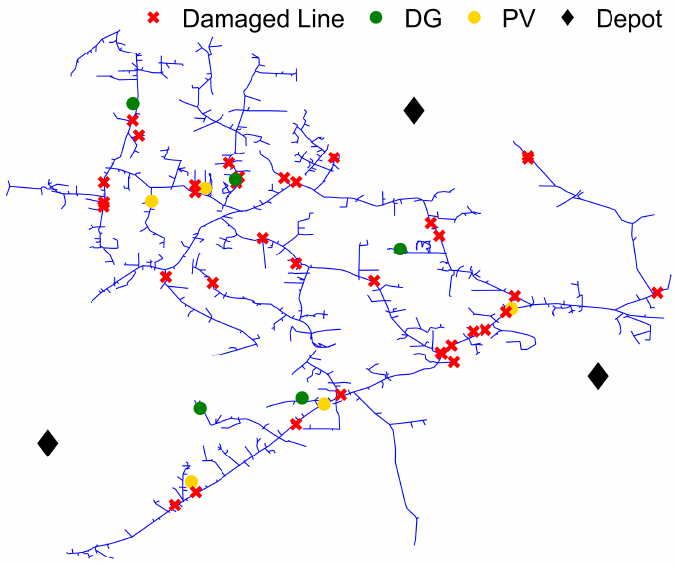}
  \caption{{\color{black}Modified IEEE 8500-bus system with 35 damaged lines.}}\label{8500}
\end{figure}

\begin{figure}[h!]
\setlength{\abovecaptionskip}{0pt} 
\setlength{\belowcaptionskip}{0pt} 
  \centering
    \includegraphics[width=0.49\textwidth]{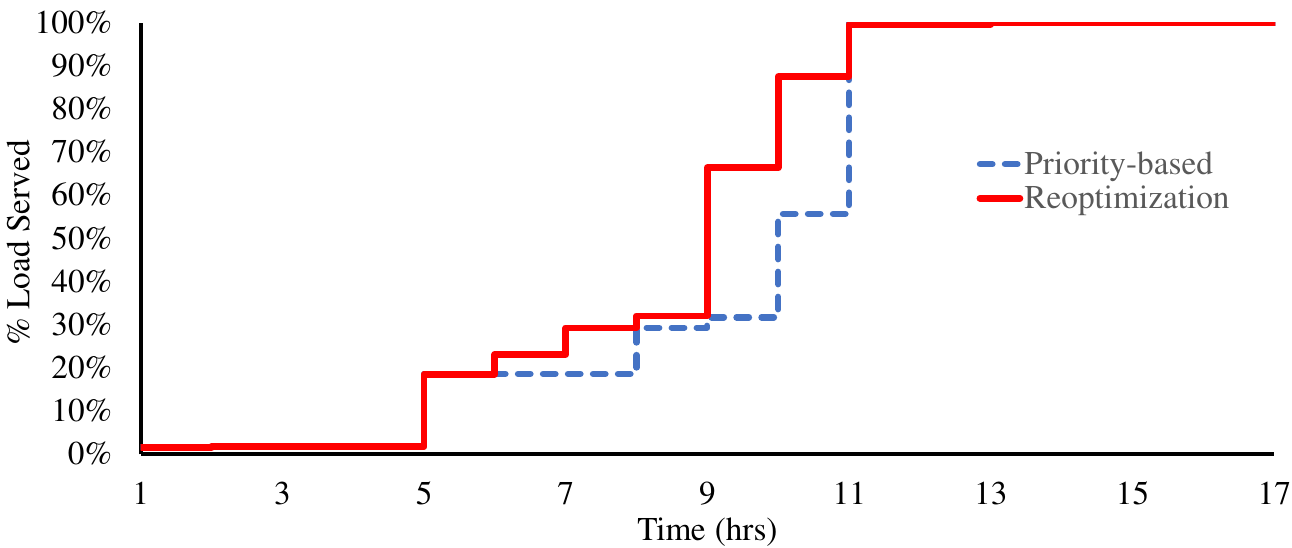}
  \caption{{\color{black}Percentage of load served at each time step for the IEEE 8500-bus system.}}\label{Load_Served_8500}
\end{figure}
}
\section{Conclusion}\label{chap:7}

In this paper, a mathematical model that combines 3-phase unbalanced distribution system operation, fault isolation and restoration, PV and BESS systems operations, and resources coordination is developed. The model included the coordination of line and tree crews as well as equipment pick up for conducting the repairs. Also, a new framework for modeling the connectivity of PV systems is designed. Furthermore, a three-stage algorithm is developed with a newly designed neighborhood search algorithm to iteratively improve the routing solution. The developed approach is able to restart when the repair time is updated, and the crews are dispatched based on the incumbent solution. Test results have shown that the proposed algorithm can provide effective restoration plans within the time limit.

%
%
\vspace{-0.8cm}
\begin{IEEEbiography}[{\includegraphics[width=1in,height=1.25in,clip,keepaspectratio]{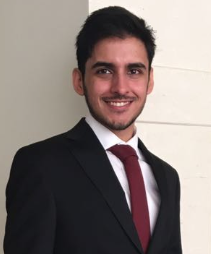}}]{Anmar Arif}
(S'16) is currently pursuing his Ph.D. in the Department of Electrical and Computer Engineering, Iowa State University, Ames, IA. He received his B.S. and Masters degrees in Electrical Engineering from King Saud University and Arizona State University in 2012 and 2015, respectively. Anmar was a Teaching Assistant in King Saud University, and a Research Assistant in Saudi Aramco Chair In Electrical Power, Riyadh, Saudi Arabia, 2013. His current research interest includes power system optimization, outage management, and microgrids.
\end{IEEEbiography}
\vspace{-0.8cm}
\begin{IEEEbiography}[{\includegraphics[width=1in,height=1.25in,clip,keepaspectratio]{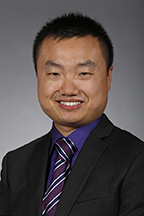}}]{Zhaoyu Wang}
(M'15) received the B.S. and M.S. degrees in electrical engineering from Shanghai Jiaotong University in 2009 and 2012, respectively, and the M.S. and Ph.D. degrees in electrical and computer engineering from the Georgia Institute of Technology in 2012 and 2015, respectively. He is the Harpole-Pentair Assistant Professor with Iowa State University. 
He was a Research Aid with Argonne National Laboratory in 2013, and an Electrical Engineer with Corning Inc. in 2014. His research interests include power distribution systems, microgrids, renewable integration, power system resilience, and power system modeling. He is the Principal Investigator for a multitude of projects focused on the above areas and funded by the National Scienc Foundation, the Department of Energy, National Laboratories, PSERC, and Iowa Energy Center. He was a recipient of the IEEE PES General Meeting Best Paper Award in 2017 and the IEEE Industrial Application Society Prize Paper Award in 2016. He is the Secretary of IEEE Power and Energy Society Award Subcommitte. He is an Editor of the IEEE TRANSACTIONS ON POWER SYSTEMS, the IEEE TRANSACTIONS ON SMART GRID and the IEEE PES LETTERS, and an Associate Editor of IET Smart Grid.
\end{IEEEbiography}
\vspace{-0.8cm}
\begin{IEEEbiography}[{\includegraphics[width=1in,height=1.25in,clip,keepaspectratio]{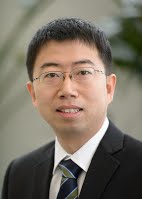}}]{Chen Chen}
(M'13) received the B.S. and M.S. degrees in electrical engineering from Xian Jiaotong University, Xian, China, in 2006 and 2009, respectively, and the Ph.D. degree in electrical engineering from Lehigh University, Bethlehem, PA, USA, in 2013. During 2013-2015, he worked as a Postdoctoral Researcher at the Energy Systems Division, Argonne National Laboratory, Argonne, IL, USA. Dr. Chen is currently an Energy Systems Scientist with the Energy Systems Division at Argonne National Laboratory. His primary research is in optimization, communications and signal processing for smart electric power systems, power system resilience, and cyber-physical system modeling for smart grids. He is an editor of IEEE Transactions on Smart Grid and IEEE Power Engineering Letters.
\end{IEEEbiography}
\vspace{-13.0cm}
\begin{IEEEbiography}[{\includegraphics[width=1in,height=1.25in,clip,keepaspectratio]{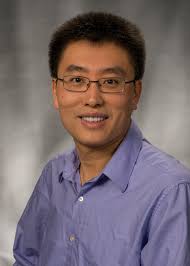}}]{Jianhui Wang}
(M'07-SM'12) received the Ph.D. degree in electrical engineering from the Illinois Institute of Technology, Chicago, IL, USA, in 2007. Presently, he is an Associate Professor with the Department of Electrical Engineering at Southern Methodist University, Dallas, Texas, USA. Prior to joining SMU, Dr. Wang had an eleven-year stint at Argonne National Laboratory with the last appointment as Section Lead - Advanced Grid Modeling. Dr. Wang is the secretary of the IEEE Power \& Energy Society (PES) Power System Operations, Planning \& Economics Committee. He has held visiting positions in Europe, Australia and Hong Kong including a VELUX Visiting Professorship at the Technical University of Denmark (DTU). Dr. Wang is the Editor-in-Chief of the IEEE Transactions on Smart Grid and an IEEE PES Distinguished Lecturer. He is also a Clarivate Analytics highly cited researcher for 2018.
\end{IEEEbiography}
\end{document}